\documentclass[11pt]{article}

\usepackage[margin=0.7in]{geometry}
\usepackage{setspace}
\usepackage{subcaption}
\usepackage[T1]{fontenc}
\usepackage[utf8]{inputenc}
\usepackage{lmodern}
\usepackage{microtype}

\usepackage{amsmath,amssymb,amsfonts,mathtools,bm, comment}
\usepackage{amsthm}
\usepackage{graphicx}
\usepackage{subcaption}
\usepackage{booktabs}
\usepackage{enumitem}
\usepackage{float}
\usepackage{algorithm}
\usepackage{algpseudocode}
\usepackage{xcolor}
\usepackage{hyperref}
\usepackage[numbers,sort&compress]{natbib}

\hypersetup{
  colorlinks=true,
  linkcolor=blue,
  citecolor=blue,
  urlcolor=blue
}

\newtheorem{theorem}{Theorem}[section]
\newtheorem{proposition}[theorem]{Proposition}

\newtheorem{corollary}[theorem]{Corollary}
\numberwithin{equation}{section}
\numberwithin{figure}{section}
\numberwithin{table}{section}

\newcommand{\R}{\mathbb{R}}
\newcommand{\dd}{\mathrm{d}}
\newcommand{\dx}{\Delta x}

\newcommand{\Eh}{\mathcal E_h}

\newcommand{\norm}[2]{\left\|#1\right\|_{#2}}
\newcommand{\inner}[2]{\left\langle #1,#2\right\rangle}

\title{\textbf{An Energy-Stable Approach for Learning Derivative Operators from Noisy Data for Maxwell’s Equations}}
\author{
Victory C. Obieke\thanks{Corresponding author. Email: \texttt{obiekev@oregonstate.edu}}\\
Oregon State University\\
\texttt{obiekev@oregonstate.edu}
\and
Ameh Emmanuel Sunday\\
Cornell University\\
\texttt{esa62@cornell.edu}
}

\date{}

\begin{document}
\maketitle

\begin{abstract}
We develop a structure-preserving ADMM method, denoted SP-ADMM, for learning
energy-stable spatial derivative stencils for Maxwell equations from noisy
data. Starting from the source-free Maxwell system, we focus on a
one-dimensional reduction whose energy conservation depends on the
skew-adjointness of the spatial derivative operator. The learned derivative is
represented by a compact periodic convolution stencil. Unlike standard
constrained ADMM, which learns the full stencil and imposes skew-adjointness
through equality constraints, SP-ADMM enforces skew-adjointness by construction
through a reduced parameterization using only the independent positive-side
stencil coefficients.

Numerical experiments show that SP-ADMM is especially effective in hidden
operator and noisy-data regimes. Across clean data, noisy derivative data,
multiple initial conditions, different hidden skew-adjoint operators,
training-set sizes, regularization parameters, constraint ablations, and
long-time simulations, SP-ADMM achieves the smallest final-time electric-field
error while preserving energy to roundoff accuracy. A layered-medium Maxwell
propagation test further shows that the learned structure-preserving stencil
remains competitive with classical finite differences in a physical
reflection/transmission setting. Overall, SP-ADMM provides a data-driven way
to learn accurate Maxwell stencils while retaining the energy-conserving
structure of the underlying equations.
\end{abstract}
\section{Introduction}
\label{sec:intro}

Finite-difference methods are a classical and reliable tool for Maxwell's
equations, including the Yee finite-difference time-domain method and standard
computational electrodynamics extensions~\cite{yee1966numerical,taflove2000computational}. Classical finite-difference stencils are fixed analytically through Taylor
moment conditions. This is appropriate when the target dynamics are generated
by the standard derivative operator. However, in data-driven settings,
model-reduction settings, or hidden-operator identification problems, the
effective spatial operator may be nonstandard. A fixed finite-difference
stencil can then be stable but inaccurate because it cannot adapt to the
observed derivative data.

This paper studies a structure-preserving alternative: learn the stencil from
data while preserving the skew-adjointness structure that gives Maxwell energy
conservation. This connects to data-driven discretization and operator
discovery~\cite{barsinai2019,brunton2016discovering}, as well as structure-preserving
scientific machine learning~\cite{obieke2025structurepreservingphysicsinformedneuralnetwork}.We begin from the source-free Maxwell curl system and then study a
one-dimensional periodic reduction that retains the skew-adjoint derivative
structure responsible for electromagnetic energy conservation. This reduced
setting provides a controlled framework for isolating the effect of
structure-preserving stencil learning before extending the approach to higher
dimensional Maxwell systems. The learned stencil is trained from derivative data generated by
hidden skew-adjoint operators. The key question is whether a learned
energy-conserving stencil can adapt to hidden spatial dynamics better than a
fixed finite-difference stencil.

We compare three methods: the fixed centered finite-difference stencil,
denoted FD; a standard constrained ADMM method that learns the full stencil
while imposing skew-adjointness as an equality constraint, denoted ADMM; and
a structure-preserving ADMM method that builds skew-adjointness into the
stencil parameterization, denoted SP-ADMM.

\begin{table}[htbp]
\centering
\caption{Qualitative comparison of fixed and learned Maxwell derivative operators.}
\label{tab:method-comparison}
\begin{tabular}{lcc}
\toprule
\textbf{Method} & \textbf{Adaptive to data?} & \textbf{Energy-preserving?} \\
\midrule
FD & \(\times\) & \(\checkmark\) \\
Unconstrained learned stencil & \(\checkmark\) & \(\times\), not guaranteed \\
ADMM & \(\checkmark\) & \(\checkmark\), by constraints \\
SP-ADMM & \(\checkmark\) & \(\checkmark\), by construction \\
\bottomrule
\end{tabular}
\end{table}

Table~\ref{tab:method-comparison} summarizes the motivation for the proposed
method. Classical finite differences provide a strong fixed baseline and
preserve the Maxwell energy structure, but they are not adaptive to data.
Unconstrained learned stencils can adapt to data, but they do not generally
preserve the skew-adjoint structure required for energy conservation. The
SP-ADMM formulation combines these two desirable features by learning from
data while remaining inside the skew-adjoint energy-conserving class by
construction.

The proposed SP-ADMM method builds on the ADMM framework~\cite{boyd2011admm}
and the broader structure-preserving viewpoint in geometric numerical
integration~\cite{lubich2006geometric}. This reformulation is not merely an implementation detail. It changes
the learned variables so that every iterate corresponds to a skew-adjoint
Maxwell operator. Thus, the energy-conserving structure is enforced by
construction.

The contributions are as follows. First, we formulate energy-conserving
Maxwell stencil learning as a constrained least-squares problem. Second, we
introduce a skew-parameterized structure-preserving ADMM solver that reduces
the optimization dimension and enforces skew-adjointness exactly. Third, we
compare FD, ADMM, and SP-ADMM over a large collection of experiments involving
noise, hidden operators, initial conditions, stencil radius, training-set size,
regularization, constraint ablation, and propagation time. We observe that SP-ADMM often improves field accuracy in the hidden-operator setting
while preserving energy to roundoff accuracy.

The results should not be interpreted as saying that learned stencils always
replace finite differences. When the data are simple, low-frequency, or
generated by a classical finite-difference operator, FD can remain
competitive. The main advantage of learning appears when the observed spatial
dynamics are generated by a \emph{nonstandard hidden skew-adjoint operator}.
Here, \emph{hidden} means that the data-generating derivative operator is unknown to the learning algorithm; and \emph{nonstandard} means that the
stencil coefficients of the derivative operator differ from those of the classical
centered finite-difference derivative. Thus, FD may still preserve energy,
but it may not match the hidden data-generating dynamics.

The rest of the paper is organized as follows.
Section~\ref{sec:stencils} introduces energy-conserving periodic convolution
stencils for the one-dimensional Maxwell system and explains why
skew-adjointness implies discrete energy conservation.
Section~\ref{sec:learning} formulates the constrained stencil-learning
problem and presents both the standard ADMM method and the proposed
structure-preserving ADMM method.
Section~\ref{sec:experiments} compares FD, ADMM, and SP-ADMM across clean
hidden-operator learning, noisy derivative data, multiple initial conditions,
different hidden operators, stencil radii, training-set sizes,training time,speed, regularization
parameters, constraint ablations,increasing final times, and a layered-medium Maxwell propagation test.
Section~\ref{sec:conclusion} concludes the paper and outlines future
extensions.

\section{Maxwell Equations and Energy-Conserving Convolution Stencils}
\label{sec:stencils}
\subsection{Continuous Maxwell system}

We begin with the source-free Maxwell equations in a homogeneous medium,
\begin{equation}
\begin{aligned}
  \varepsilon \frac{\partial \mathbf E}{\partial t}
  &= \nabla \times \mathbf H,\\
  \mu \frac{\partial \mathbf H}{\partial t}
  &= - \nabla \times \mathbf E,
\end{aligned}
\qquad \mathbf x\in \Omega,\quad t>0,
\label{eq:full-maxwell}
\end{equation}
where $\mathbf E(\mathbf x,t)$ is the electric field, $\mathbf H(\mathbf x,t)$
is the magnetic field, and $\varepsilon,\mu>0$ are the electric permittivity
and magnetic permeability. For periodic boundary conditions, the continuous
electromagnetic energy
\begin{equation}
  \mathcal E(t)
  =
  \frac12 \int_{\Omega}
  \left(
  \varepsilon |\mathbf E(\mathbf x,t)|^2
  +
  \mu |\mathbf H(\mathbf x,t)|^2
  \right)\,d\mathbf x
  \label{eq:continuous-energy}
\end{equation}
is conserved. Indeed, using integration by parts and the skew-adjointness of
the curl operator under periodic boundary conditions,
\[
  \frac{d\mathcal E}{dt}
  =
  \int_{\Omega}
  \mathbf E\cdot(\nabla\times \mathbf H)
  -
  \mathbf H\cdot(\nabla\times \mathbf E)
  \,d\mathbf x
  =
  0.
\]

The one-dimensional model used in the present stencil-learning experiments is
obtained by considering fields depending only on $x$ and choosing a transverse
electric and magnetic component. After nondimensionalization, this gives
\begin{equation}
  E_t = H_x,\qquad H_t = E_x,
  \qquad x\in[0,L],\quad t>0.
  \label{eq:oned-maxwell-model}
\end{equation}
Thus, the one-dimensional system studied below is a reduced Maxwell system
that retains the essential skew-adjoint derivative structure responsible for
energy conservation. The learning problem considered in this work is to learn
a discrete derivative operator that preserves this structure.

\subsection{Periodic convolution derivative operators}

Let $x_i=i\dx$, $i=0,\ldots,N-1$, be a periodic grid on $[0,L]$, with
$\dx=L/N$. A periodic convolution derivative stencil of radius $R$ is
\begin{equation}
  (D_w u)_i = \sum_{\ell=-R}^{R} w_\ell u_{i+\ell},
  \qquad i=0,\ldots,N-1,
  \label{eq:stencil}
\end{equation}
The coefficient vector is
\[
  w = (w_{-R},\ldots,w_{-1},w_0,w_1,\ldots,w_R)^T.
\]

For periodic convolution matrices, skew-adjointness is equivalent to
antisymmetry of the stencil:
\begin{equation}
  w_0=0,\qquad w_{-\ell}=-w_\ell,\qquad \ell=1,\ldots,R.
  \label{eq:skew}
\end{equation}
When \eqref{eq:skew} holds, $D_w^T=-D_w$.

The semi-discrete Maxwell system associated with $D_w$ is
\begin{equation}
  \frac{\dd E}{\dd t}=D_wH,
  \qquad
  \frac{\dd H}{\dd t}=D_wE.
  \label{eq:semidiscrete}
\end{equation}
We define the discrete electromagnetic energy by
\begin{equation}
  \Eh(t)=
  \frac12\norm{E(t)}{\ell^2_h}^2
  +
  \frac12\norm{H(t)}{\ell^2_h}^2,
  \qquad
  \norm{u}{\ell^2_h}^2=\dx\sum_i u_i^2.
  \label{eq:energy}
\end{equation}
We recall the standard discrete energy conservation property associated
with skew-adjoint derivative operators, \cite{bokil2014operator,obieke2025structurepreservingphysicsinformedneuralnetwork}.
\begin{proposition}[Discrete energy conservation]
If $w$ satisfies \eqref{eq:skew}, then the semi-discrete Maxwell system
\eqref{eq:semidiscrete} conserves the energy \eqref{eq:energy}.
\end{proposition}

\begin{proof}
Differentiating \eqref{eq:energy} gives
\[
  \frac{\dd \Eh}{\dd t}
  =
  \inner{D_wH}{E}_{\ell^2_h}
  +
  \inner{D_wE}{H}_{\ell^2_h}.
\]
Since $D_w^T=-D_w$,
\[
  \inner{D_wH}{E}_{\ell^2_h}
  =
  -\inner{H}{D_wE}_{\ell^2_h}
  =
  -\inner{D_wE}{H}_{\ell^2_h}.
\]
The two terms cancel, so $\dd\Eh/\dd t=0$.
\end{proof}

This proposition is the foundation of the learning framework. The learned
operator may be data-adaptive, but as long as it is skew-adjoint, it preserves
the Maxwell energy structure.
\begin{corollary}[Energy stability]
If $w$ satisfies \eqref{eq:skew}, then the semi-discrete Maxwell solution is
uniformly bounded in the discrete energy norm. In particular,
\[
  \norm{E(t)}{\ell^2_h}^2+\norm{H(t)}{\ell^2_h}^2
  =
  \norm{E(0)}{\ell^2_h}^2+\norm{H(0)}{\ell^2_h}^2,
  \qquad t\ge 0.
\]
\end{corollary}

\begin{proof}
This follows directly from the energy conservation identity
$\Eh(t)=\Eh(0)$ and the definition of $\Eh$ in \eqref{eq:energy}.
\end{proof}

\section{Constrained Stencil Learning}
\label{sec:learning}

\subsection{Training system}

Training data are built from sampled states $(E^{(m)},H^{(m)})$ and target
time derivatives. For Maxwell's equations, the targets have the form
\[
  b_E^{(m)}\approx \partial_t E^{(m)}=D_{\rm true}H^{(m)},\qquad
  b_H^{(m)}\approx \partial_t H^{(m)}=D_{\rm true}E^{(m)}.
\]
Stacking local stencil equations over grid points and training samples gives
a linear system
\begin{equation}
  Aw \approx b.
  \label{eq:Awb}
\end{equation}
In the hidden-operator experiments, $b$ is generated by a nonstandard
skew-adjoint operator $D_{\rm true}$, and then Gaussian noise is added:
\begin{equation}
  b_{\rm noisy}=b_{\rm clean}+\eta,\qquad
  \eta\sim \mathcal N(0,\sigma^2I),\qquad
  \sigma=\delta\,{\rm std}(b_{\rm clean}),
  \label{eq:noise}
\end{equation}
where $\delta$ is the noise level.

\subsection{Standard constrained ADMM}

The standard constrained learning problem is
\begin{equation}
\begin{aligned}
  \min_{w\in\R^{2R+1}}
  \quad&
  \frac12\norm{Aw-b}{2}^2+\frac{\lambda}{2}\norm{w}{2}^2,\\
  \text{subject to}\quad&
  Cw=0,\\
  & -M_b\le w_\ell\le M_b.
\end{aligned}
\label{eq:full-qP}
\end{equation}
Here $Cw=0$ represents the skew-adjointness constraints
\[
  w_0=0,\qquad w_{-\ell}+w_\ell=0.
\]
Additional moment constraints can be included when the goal is to recover a
classical finite-difference stencil.

A standard ADMM splitting introduces an auxiliary variable $z$ for the box
constraint and solves equality-constrained quadratic subproblems for $w$.
The $w$-update is obtained from the KKT system
\begin{equation}
\begin{bmatrix}
A^TA+(\lambda+\rho)I & G^T\\
G & 0
\end{bmatrix}
\begin{bmatrix}
w^{k+1}\\
\mu^{k+1}
\end{bmatrix}
=
\begin{bmatrix}
A^Tb+\rho(z^k-u^k)\\
h
\end{bmatrix},
\label{eq:kkt}
\end{equation}
where $Gw=h$ contains the active equality constraints.

\subsection{Structure-preserving ADMM}
As illustrated in Figure 2.1, the proposed SP-ADMM method uses the fact that skew-adjointness can be built directly into the stencil parameterization. Let
\[
  a=(a_1,\ldots,a_R)^T
\]
contain only the positive-side stencil coefficients, and define a matrix
$S\in\R^{(2R+1)\times R}$ such that
\begin{equation}
  w=Sa=
  (-a_R,\ldots,-a_2,-a_1,0,a_1,a_2,\ldots,a_R)^T.
  \label{eq:skew-param}
\end{equation}
Then $w$ is skew-symmetric for every $a$. The constrained problem becomes
\begin{equation}
\begin{aligned}
  \min_{a\in\R^R}
  \quad&
  \frac12\norm{ASa-b}{2}^2+\frac{\lambda}{2}\norm{a}{2}^2,\\
  \text{subject to}\quad&
  -M_b\le a_j\le M_b.
\end{aligned}
\label{eq:sp-problem}
\end{equation}
The skew constraint no longer appears as an equality constraint; it is
satisfied by construction.

Introducing an auxiliary variable $z$ for the box constraint gives
\begin{align}
  a^{k+1}
  &=
  \left((AS)^T(AS)+(\lambda+\rho)I\right)^{-1}
  \left((AS)^Tb+\rho(z^k-u^k)\right),
  \label{eq:sp-a-update}\\
  \widehat a^{k+1}
  &=
  \alpha a^{k+1}+(1-\alpha)z^k,
  \label{eq:sp-relax}\\
  z^{k+1}
  &=
  \Pi_{[-M_b,M_b]^R}(\widehat a^{k+1}+u^k),
  \label{eq:sp-z-update}\\
  u^{k+1}
  &=
  u^k+\widehat a^{k+1}-z^{k+1}.
  \label{eq:sp-u-update}
\end{align}
Here $\alpha\in(1,2)$ is an over-relaxation parameter. In the experiments,
an adaptive penalty parameter is used to balance primal and dual residuals.

Algorithm~\ref{alg:sp-admm} summarizes the skew-parameterized
structure-preserving ADMM iteration.

\begin{algorithm}
\caption{Skew-parameterized structure-preserving ADMM}
\label{alg:sp-admm}
\begin{algorithmic}[1]
\Require Design matrix $A$, target vector $b$, skew map $S$, regularization
$\lambda$, box bound $M_b$, penalty $\rho$, relaxation $\alpha$.
\State Define $A_S=AS$.
\State Initialize $z^0,u^0\in\R^R$.
\For{$k=0,1,\ldots,K-1$}
  \State Compute $a^{k+1}$ from \eqref{eq:sp-a-update}.
  \State Over-relax using \eqref{eq:sp-relax}.
  \State Project onto the box using \eqref{eq:sp-z-update}.
  \State Update the scaled dual variable using \eqref{eq:sp-u-update}.
  \State Optionally update $\rho$ using residual balancing.
\EndFor
\State Return the skew stencil $w_{\rm SP}=Sz^K$.
\end{algorithmic}
\end{algorithm}
Since every SP-ADMM stencil satisfies the skew-adjointness condition
\eqref{eq:skew} by construction, the corresponding semi-discrete Maxwell
system is energy-stable by the corollary.
The advantage of SP-ADMM is structural: every iterate remains inside the energy-conserving class of skew-adjoint convolution operators. ADMM learns the full stencil and
imposes skew-adjointness as a constraint. SP-ADMM only learns the independent
positive-side coefficients, so every iterate produces a skew-adjoint
convolution matrix.

\section{Numerical Experiments}
\label{sec:experiments}

We consider the nondimensional first-order Maxwell system
\begin{equation}
  E_t = H_x, \qquad H_t = E_x,
  \qquad x\in[0,L], \quad t>0,
  \label{eq:oned-maxwell}
\end{equation}
with periodic boundary conditions
\begin{equation}
  E(x+L,t)=E(x,t), \qquad H(x+L,t)=H(x,t).
  \label{eq:periodic-bc}
\end{equation}
The reference
solution is generated with a hidden skew-adjoint operator.
Numerical solutions are generated using FD, ADMM, and SP-ADMM operators. To focus on spatial discretization error, the semi-discrete systems are
advanced using the matrix exponential, a standard tool for the exact
evolution of linear constant-coefficient systems and for numerical
matrix-function computations \cite{higham2008functions,moler2003dubious}. We measure the final-time electric-field error and relative energy drift by
\[
  {\rm err}_E(T)=
  \frac{\norm{E_M(T)-E_{\rm true}(T)}{\ell^2_h}}
       {\norm{E_{\rm true}(T)}{\ell^2_h}},
  \qquad
  {\rm drift}(T)=
  \frac{|\Eh(T)-\Eh(0)|}{|\Eh(0)|}.
\]
Values below $10^{-18}$ are plotted at $10^{-18}$ for visualization on
logarithmic axes.

\subsection{Finite-difference recovery}
\label{subsec:fd-recovery}

\begin{figure}
\centering
\includegraphics[width=0.78\textwidth]{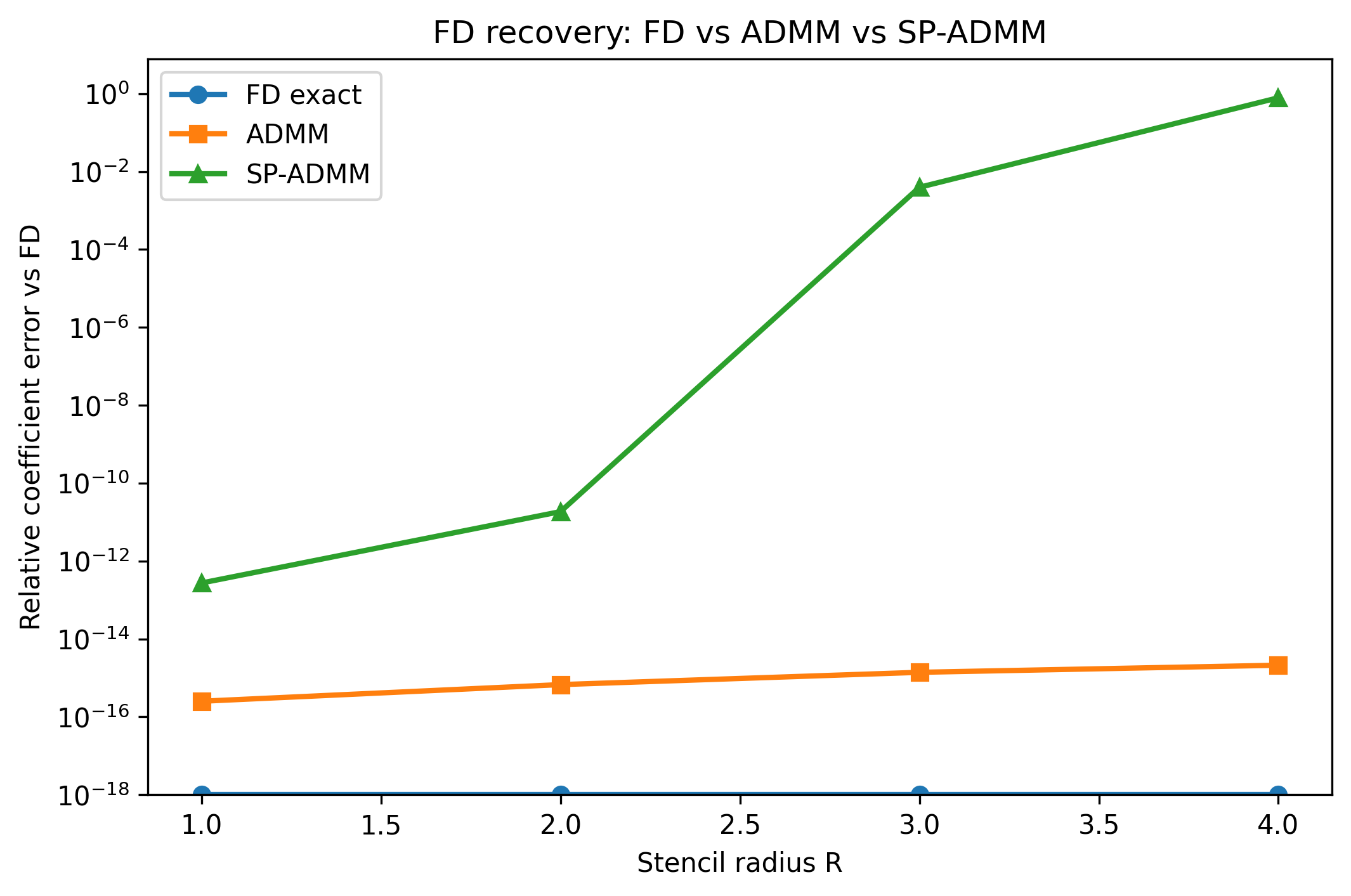}
\caption{Finite-difference recovery test. ADMM with moment constraints
recovers the classical finite-difference stencil to roundoff accuracy. The
SP-ADMM curve is a skew-parameterized learned stencil without the same moment
constraints; its purpose here is comparison, not exact FD recovery.}
\label{fig:fd-recovery}
\end{figure}

The first experiment is a verification test. The reference operator is the
classical finite-difference stencil. The ADMM method is given the appropriate
moment constraints and recovers the finite-difference coefficients to
roundoff accuracy.The SP-ADMM curve is included only for comparison. In this experiment, SP-ADMM is not supplied with the finite-difference moment constraints and therefore is not expected to recover the classical finite-difference coefficients.
Figure~\ref{fig:fd-recovery} shows the finite-difference recovery test. The
results confirm that the constrained ADMM formulation exactly reproduces the
standard finite-difference stencil up to machine precision when the finite-difference moment conditions
are enforced. Without these moment constraints, the learned stencil
is not restricted to coincide with the classical finite-difference coefficients.
This flexibility is important because it allows the proposed optimization
framework to learn structure-preserving operators that may differ from standard
finite differences while still satisfying the desired skew-adjointness and
energy-conservation constraints.

\subsection{Clean hidden-operator learning}
\label{subsec:clean}

\begin{figure}
\centering
\includegraphics[width=0.78\textwidth]{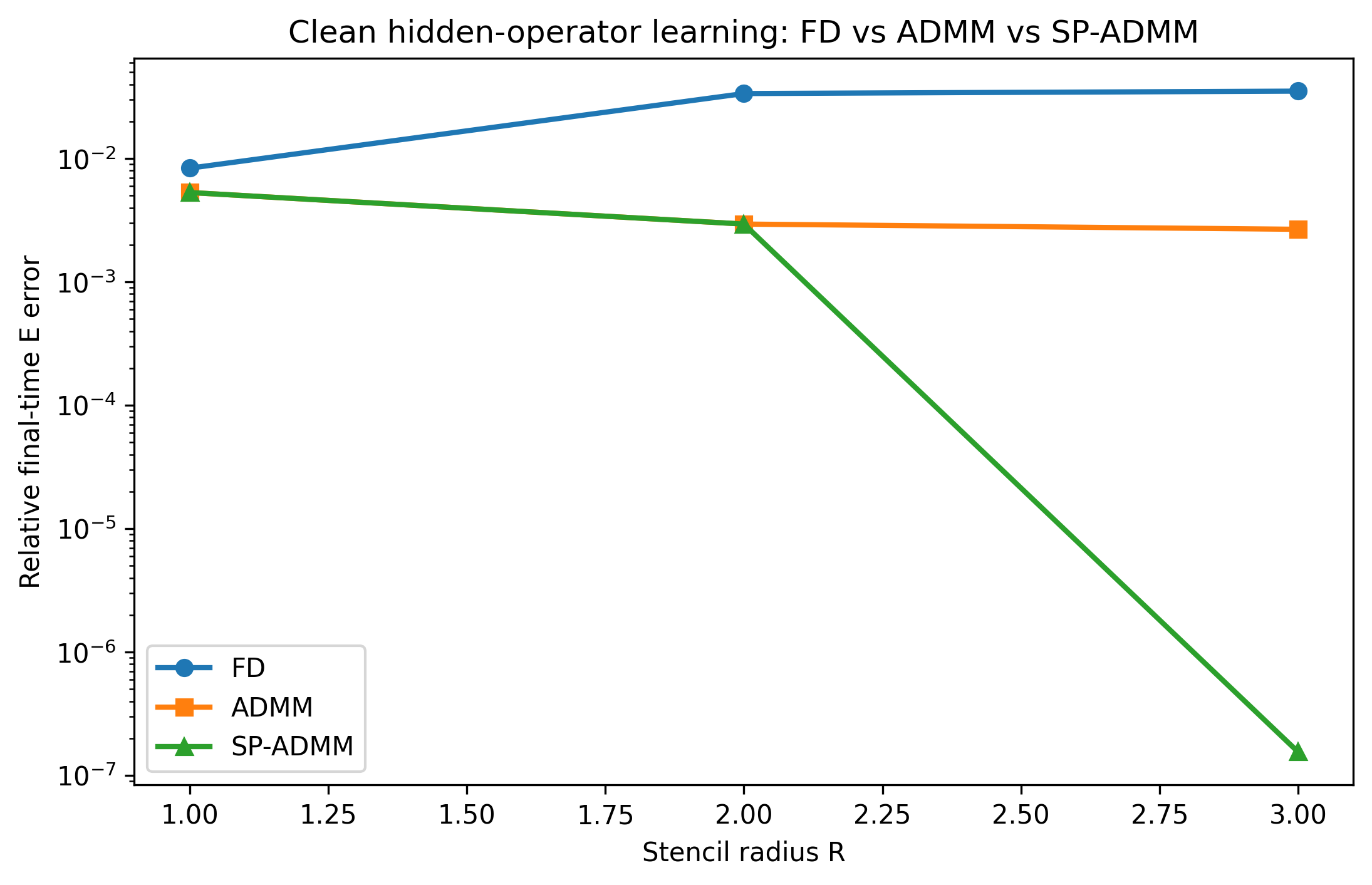}
\caption{Clean hidden-operator learning. The learned stencils outperform
fixed finite differences, and SP-ADMM gives the smallest error for several
radii. This indicates that the skew-parameterized solver can adapt strongly
to the hidden skew-adjoint operator.}
\label{fig:clean-hidden}
\end{figure}

We next train on clean derivative data generated by a hidden nonstandard
skew-adjoint operator. The learned methods are compared with fixed finite
differences for several stencil radii. Figure~\ref{fig:clean-hidden}
compares FD, ADMM, and SP-ADMM for clean hidden-operator data across
different stencil radii.

For R=1, FD, ADMM, and SP-ADMM have comparable errors. For R=2, both learned methods reduce the error relative to FD. For R=3, SP-ADMM gives a substantially smaller error, suggesting that this stencil width is sufficient to capture the hidden operator in this experiment. Therefore, unless otherwise stated, we use R=3 in the remaining experiments

\subsection{Noise robustness}
\label{subsec:noise}

\begin{figure}[htbp]
\centering

\begin{subfigure}[t]{0.32\textwidth}
    \centering
    \includegraphics[width=\textwidth]{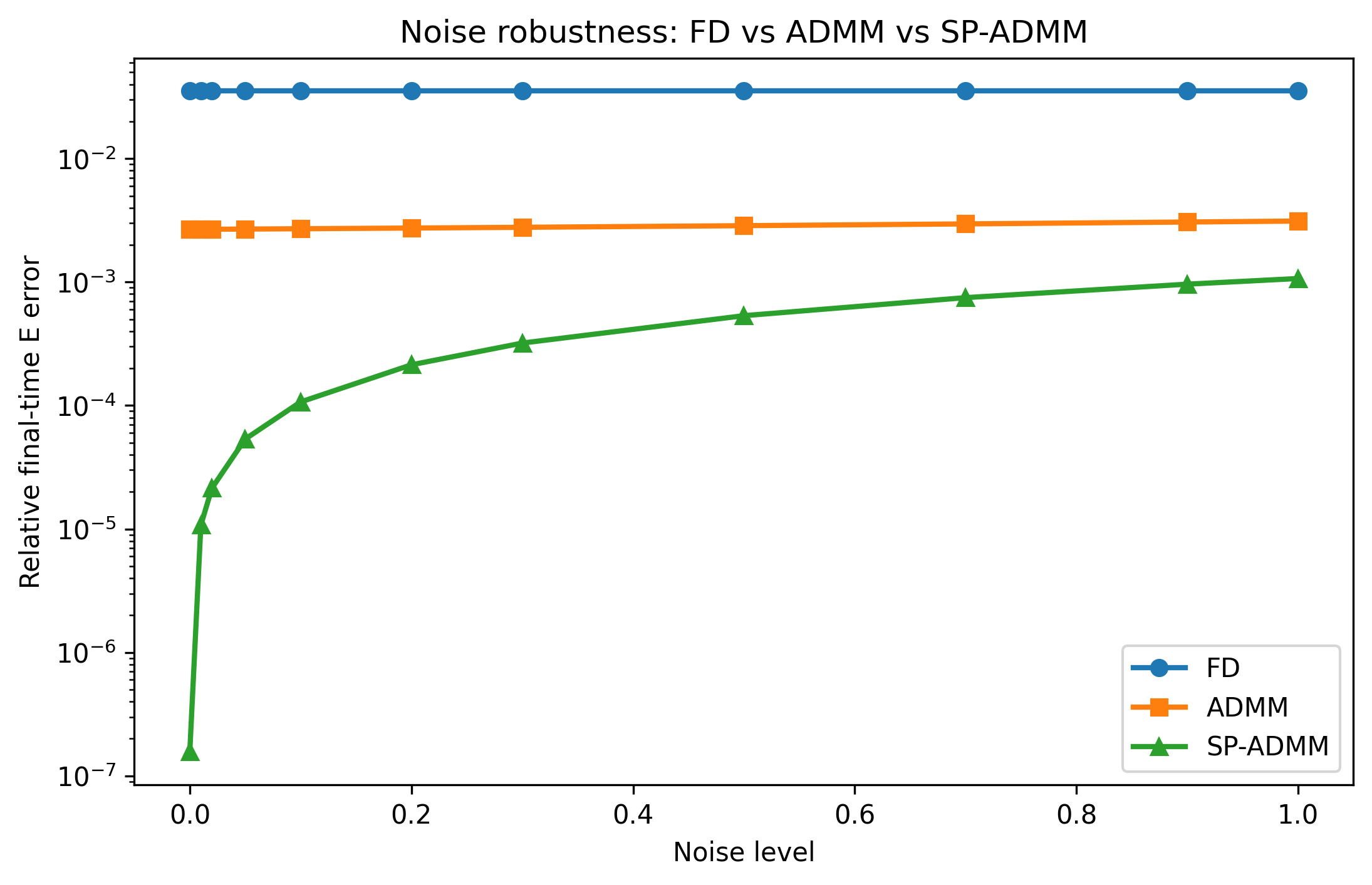}
    \caption{Final-time electric-field error.}
    \label{fig:noise-error}
\end{subfigure}
\hfill
\begin{subfigure}[t]{0.32\textwidth}
    \centering
    \includegraphics[width=\textwidth]{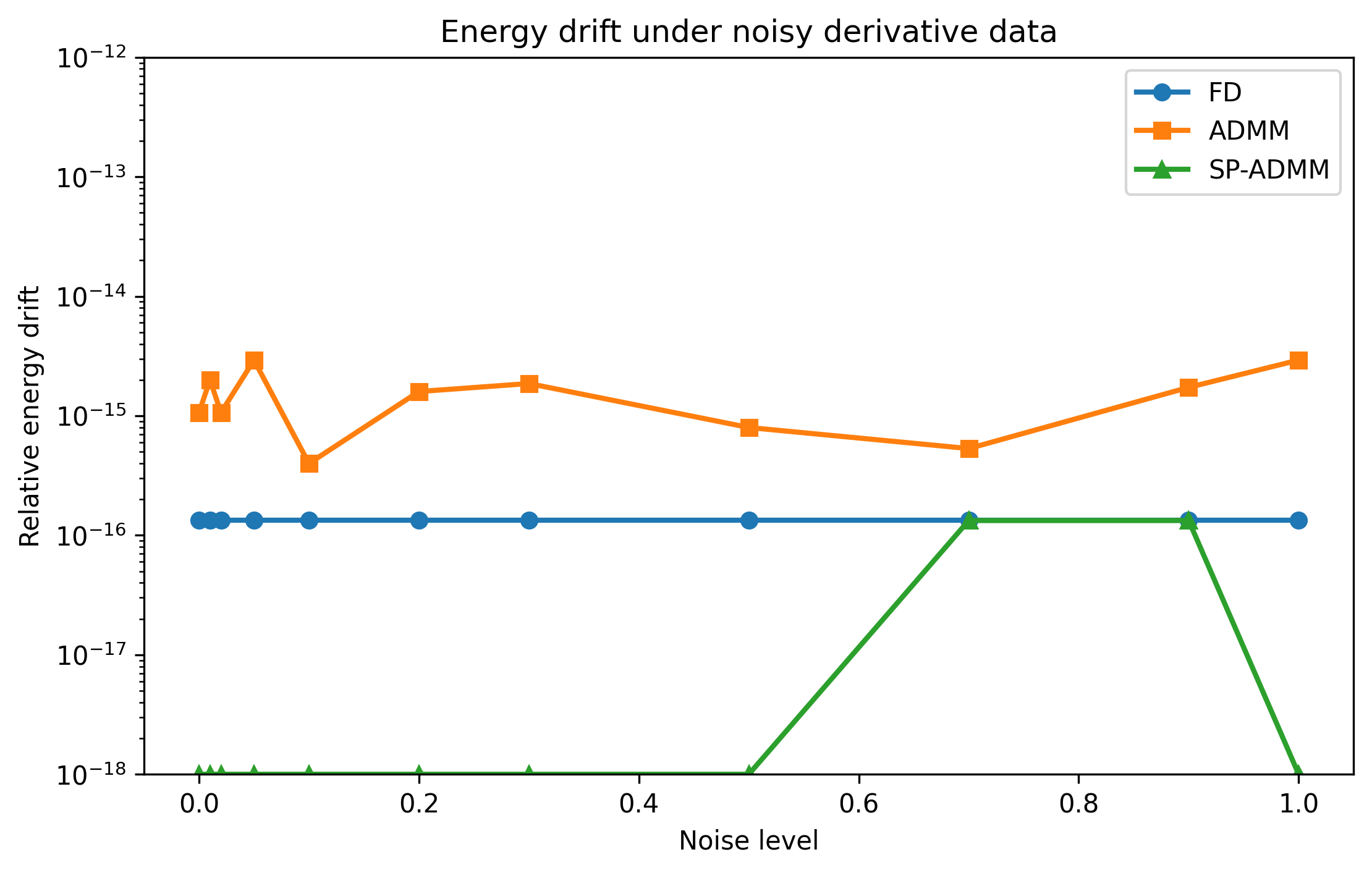}
    \caption{Energy drift.}
    \label{fig:noise-energy}
\end{subfigure}
\hfill
\begin{subfigure}[t]{0.32\textwidth}
    \centering
    \includegraphics[width=\textwidth]{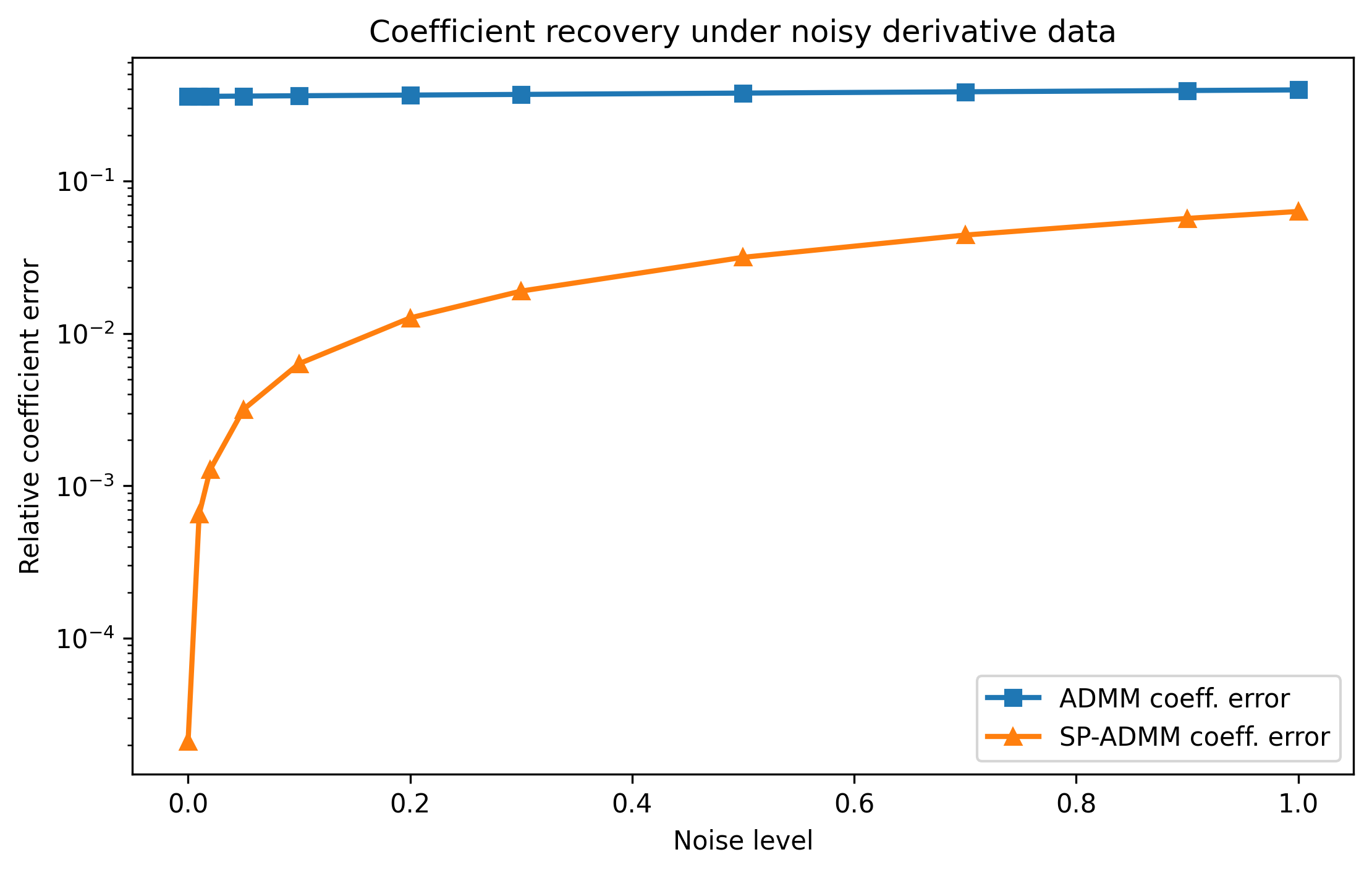}
    \caption{Coefficient recovery.}
    \label{fig:noise-coeff}
\end{subfigure}

\caption{Noise robustness under derivative-data corruption. 
(a) FD remains flat because its coefficients are fixed, while both learned methods remain below FD across the tested noise levels; SP-ADMM gives the smallest final-time electric-field error across the full noise range. 
(b) All skew-adjoint methods preserve the discrete Maxwell energy to roundoff accuracy, and the small variations between curves are due to numerical roundoff. 
(c) SP-ADMM gives smaller coefficient error than ADMM at low and moderate noise levels. As the noise level increases, both learned methods degrade, but the skew-parameterized SP-ADMM formulation remains more accurate over most of the noise range.}
\label{fig:noise-all}
\end{figure}

The noise sweep tests derivative-data corruption over many noise levels.
Finite differences are fixed and therefore independent of the training noise,
while ADMM and SP-ADMM are retrained for each noise level.
Figure~\ref{fig:noise-error} shows the final-time electric-field error,
Figure~\ref{fig:noise-energy} reports the corresponding relative energy drift,
and Figure~\ref{fig:noise-coeff} compares the coefficient-recovery error for
ADMM and SP-ADMM.

Figure~\ref{fig:noise-error} shows that SP-ADMM gives the smallest
final-time electric-field error across all tested noise levels. Since FD is
fixed, its error is independent of the training noise and remains nearly
constant. ADMM also improves substantially over FD, but SP-ADMM is more
accurate because the skew-adjoint structure is built directly into the stencil
parameterization. As the noise level increases, the SP-ADMM error increases
gradually, which is expected because the derivative data become less reliable.
However, it remains below both FD and ADMM over the full noise range.

Figure~\ref{fig:noise-energy} shows that SP-ADMM preserves the discrete
Maxwell energy to roundoff accuracy for all tested noise levels. For most
noise levels, the SP-ADMM energy drift lies at the plotting floor, indicating
that the measured drift is below the numerical resolution used for the
logarithmic plot. At a few intermediate noise levels, the SP-ADMM curve rises
to about $10^{-16}$, which is still at machine-precision level and comparable
to the FD energy drift. Therefore, the small oscillations in the SP-ADMM
energy curve should be interpreted as roundoff effects, not as physical
energy growth or instability. Compared with ADMM, which remains around
$10^{-15}$, SP-ADMM gives slightly smaller energy drift overall while also
maintaining the lowest field error in Figure~\ref{fig:noise-error}.

The coefficient-recovery diagnostic in Figure~\ref{fig:noise-coeff} shows a
similar trend. SP-ADMM has much smaller coefficient error than ADMM at low
and moderate noise levels. The coefficient error increases as the noise level
grows, but Figure~\ref{fig:noise-error} shows that SP-ADMM still produces the
most accurate Maxwell evolution. This suggests that the skew-parameterized
formulation acts as a useful structure-preserving regularization: it does not
merely fit the noisy derivative data, but learns a stable energy-conserving
stencil that gives accurate field propagation.
\subsection{Generalization across initial conditions}
\label{subsec:ics}

\begin{figure}
\centering
\includegraphics[width=0.95\textwidth]{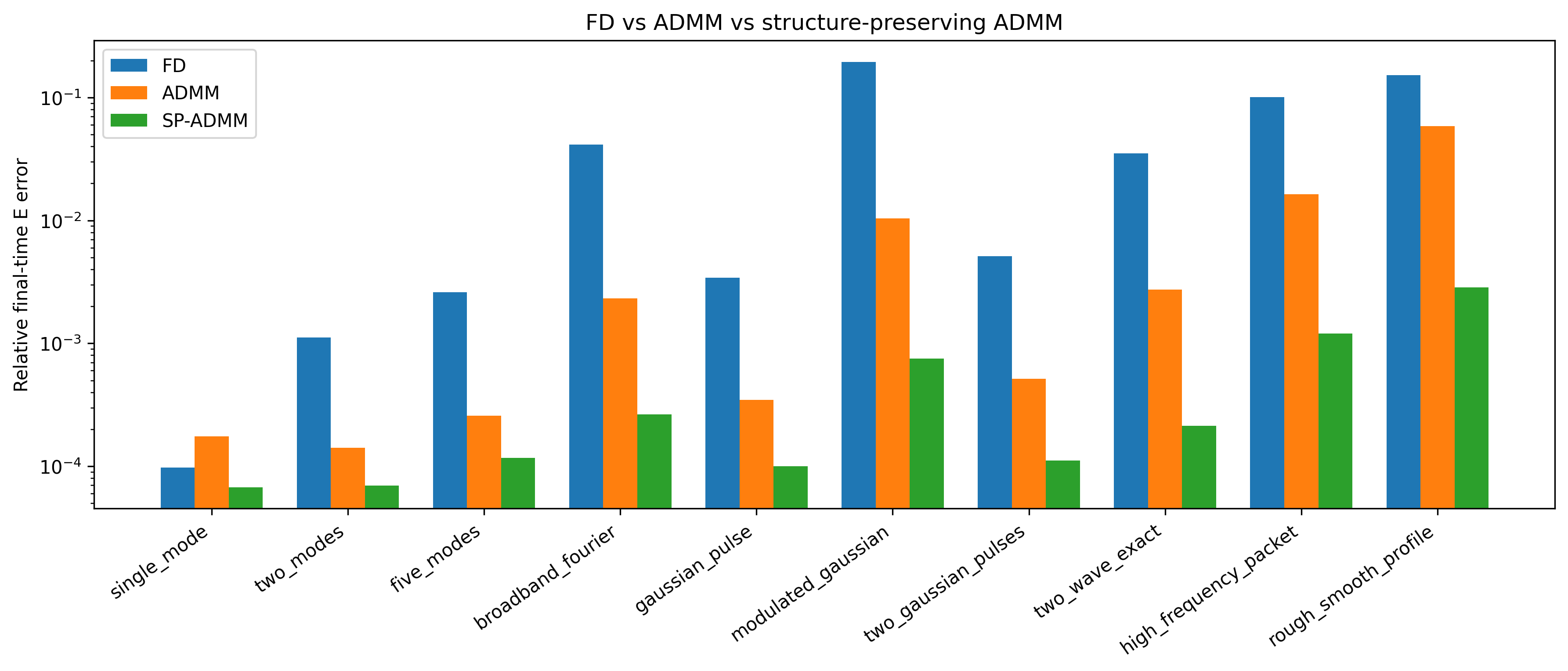}
\caption{Generalization across initial conditions. The learned ADMM and
SP-ADMM stencils reduce the error relative to the fixed finite-difference
stencil for most test profiles. SP-ADMM gives the smallest error for several
localized and broadband initial conditions, showing the benefit of enforcing
the skew-adjoint structure through the stencil parameterization.}
\label{fig:initial-conditions}
\end{figure}

We next test the learned stencils on several initial conditions, including
single-mode waves, multi-mode Fourier profiles, Gaussian pulses, localized
wave packets, and rough smooth profiles. In this experiment, the ADMM and
SP-ADMM stencils are trained on the same noisy derivative data, and both are
compared against the fixed finite-difference stencil. Figure~\ref{fig:initial-conditions}
shows the resulting final-time electric-field errors.

The results show that the learned stencils generalize beyond the training
profiles. The improvement is most pronounced for broadband, localized, and
high-frequency initial conditions, where the fixed finite-difference stencil
does not match the hidden data-generating operator as well. For very simple
low-frequency data, finite differences may remain competitive, which is
expected because classical centered stencils are already accurate for smooth
single-mode waves. Nevertheless, SP-ADMM gives consistently strong
performance across all tested initial conditions and achieves the smallest
error in all cases.

\subsection{Performance across many hidden operators}
\label{subsec:hidden-operators}

\begin{figure}
\centering
\includegraphics[width=0.90\textwidth]{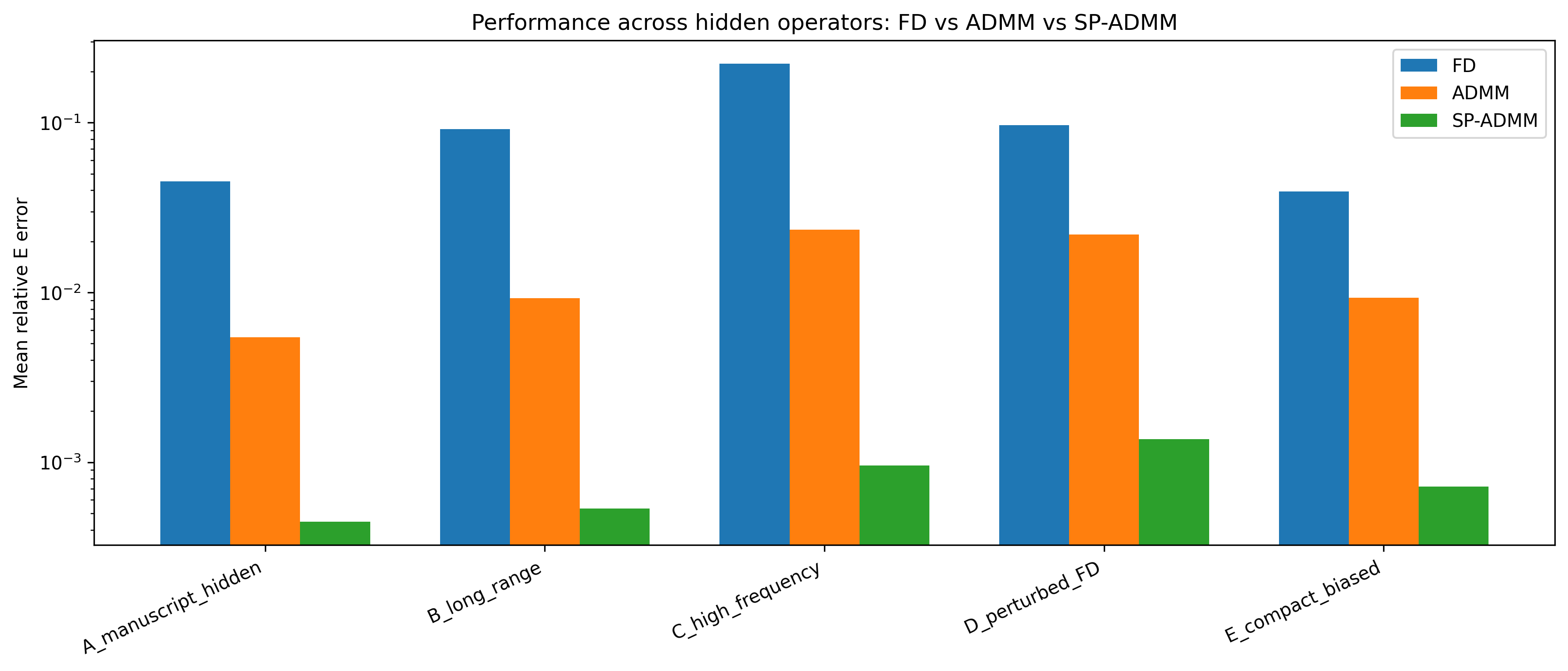}
\caption{Performance across hidden skew-adjoint operators. SP-ADMM gives the
lowest mean error for all hidden operators tested, while ADMM also improves
substantially over fixed finite differences. This shows that the advantage of
the learned stencils is not restricted to one particular hidden operator.}
\label{fig:hidden-operators}
\end{figure}

To assess robustness, we test whether the observed improvement is tied to a single data-generating operator. We generate training data from five different hidden skew-adjoint operators. Figure~\ref{fig:hidden-operators} shows the mean final-time
electric-field error for FD, ADMM, and SP-ADMM across these hidden operators.

Figure~\ref{fig:hidden-operators} shows that the fixed finite-difference
stencil has the largest error for every hidden operator. This is expected
because FD is energy-conserving but is not adapted to the nonstandard
data-generating stencil. ADMM reduces the error substantially by learning from
the derivative data, but SP-ADMM gives the smallest mean error in all cases.
The improvement is especially pronounced for the long-range, high-frequency,
and compact-biased hidden operators, where matching the effective discrete
dispersion is important.

This experiment is one of the strongest pieces of evidence for the
structure-preserving learned-stencil approach. The results show that the observed improvement is not tied to one particular hidden operator. Rather, the skew-parameterized SP-ADMM formulation consistently learns accurate energy-conserving stencils across several hidden skew-adjoint dynamics.

\subsection{Effect of training-set size}
\label{subsec:training-size}

\begin{figure}
\centering
\includegraphics[width=0.78\textwidth]{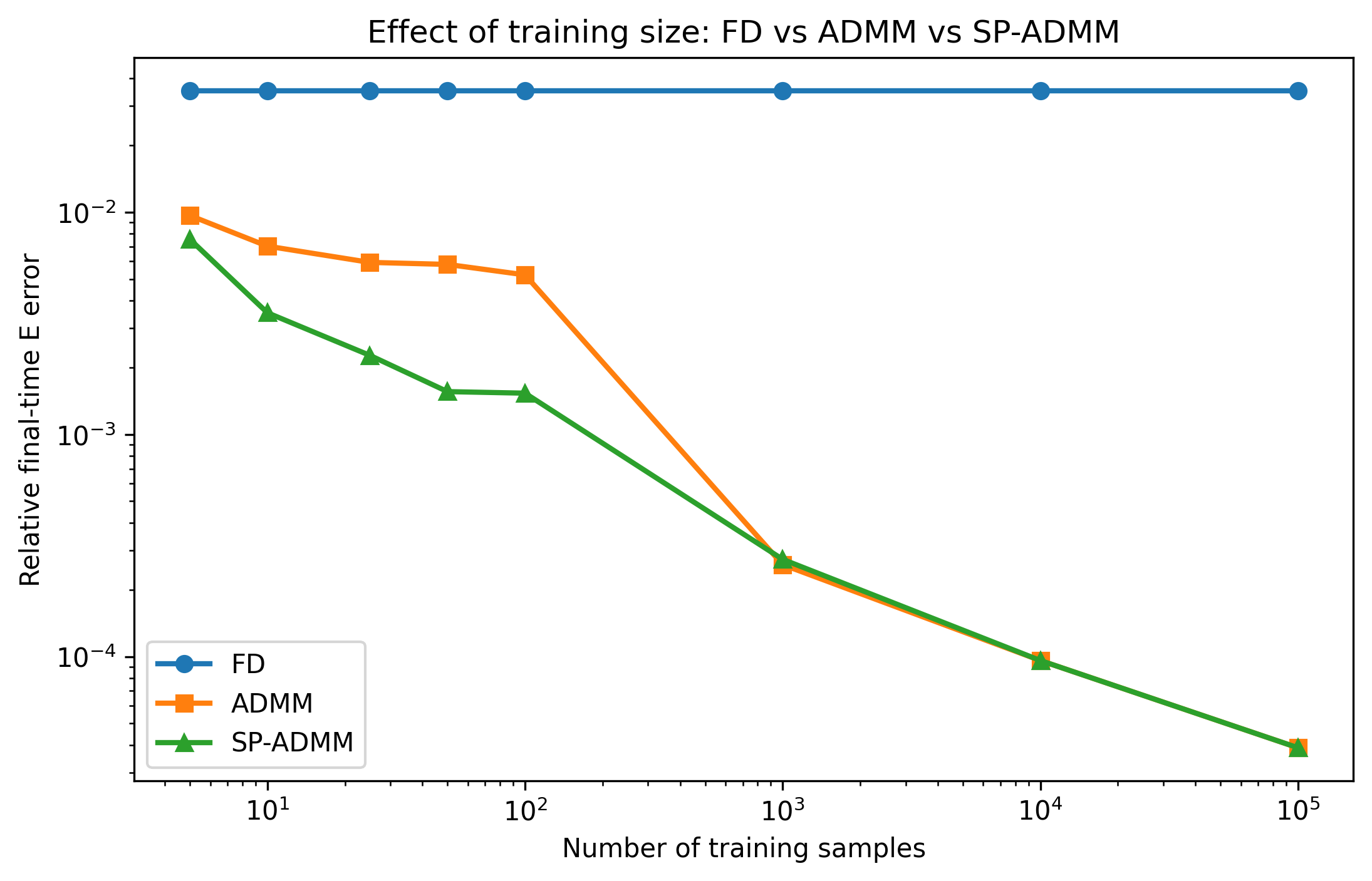}
\caption{Effect of training-set size. The FD error remains fixed because the
finite-difference stencil is not trained. The ADMM and SP-ADMM errors decrease
as more training data are added. SP-ADMM gives the smaller error in the
small- and moderate-data regimes, while both learned methods become nearly
indistinguishable for very large training sets.}
\label{fig:training-size}
\end{figure}

We next vary the number of training samples used to learn the ADMM and
SP-ADMM stencils. This experiment tests whether the learned stencils improve
as more derivative data become available. Figure~\ref{fig:training-size}
shows the final-time electric-field error as the number of training samples
increases from a small data regime to a large data regime.

Figure~\ref{fig:training-size} shows that the fixed finite-difference error is
nearly constant, as expected, because FD does not depend on the training set.
In contrast, both learned methods improve substantially as the number of
training samples increases. In the small-data regime, SP-ADMM gives a clear
advantage over ADMM, indicating that the skew-parameterized formulation
provides a useful structural bias when the available data are limited. This is
important because the admissible stencil is constrained to be skew-adjoint by
construction, reducing the effective search space and preventing the learned
operator from fitting nonphysical components of the data.

As the number of samples increases, the errors of both learned methods
continue to decrease by nearly two orders of magnitude. For large training
sets, ADMM and SP-ADMM become very close, suggesting that with enough data the
constrained ADMM formulation can also identify a highly accurate
energy-conserving stencil with SP-ADMM offering gain in training time for larger training set. Nevertheless, SP-ADMM reaches low error earlier and
remains competitive throughout the sweep. Overall, this experiment shows that
the learned stencil approach benefits strongly from additional derivative
data, while the fixed FD stencil cannot improve because its coefficients are
prescribed in advance.

\subsection{Training-time comparison}
\label{subsec:training-time}

\begin{figure}[htbp]
\centering

\begin{subfigure}[t]{0.32\textwidth}
    \centering
    \includegraphics[width=\textwidth]{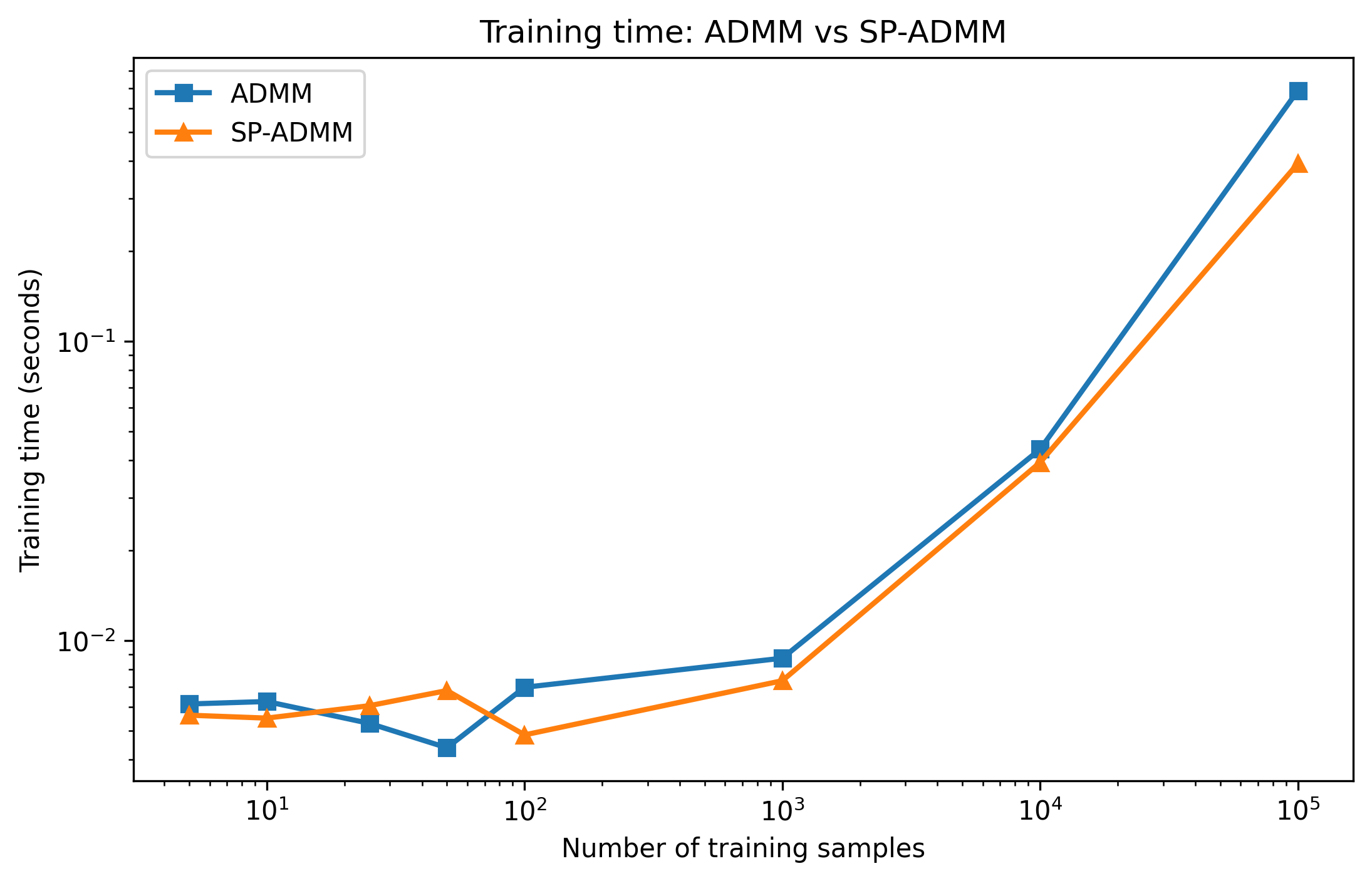}
    \caption{Total training time.}
    \label{fig:training-time-a}
\end{subfigure}
\hfill
\begin{subfigure}[t]{0.32\textwidth}
    \centering
    \includegraphics[width=\textwidth]{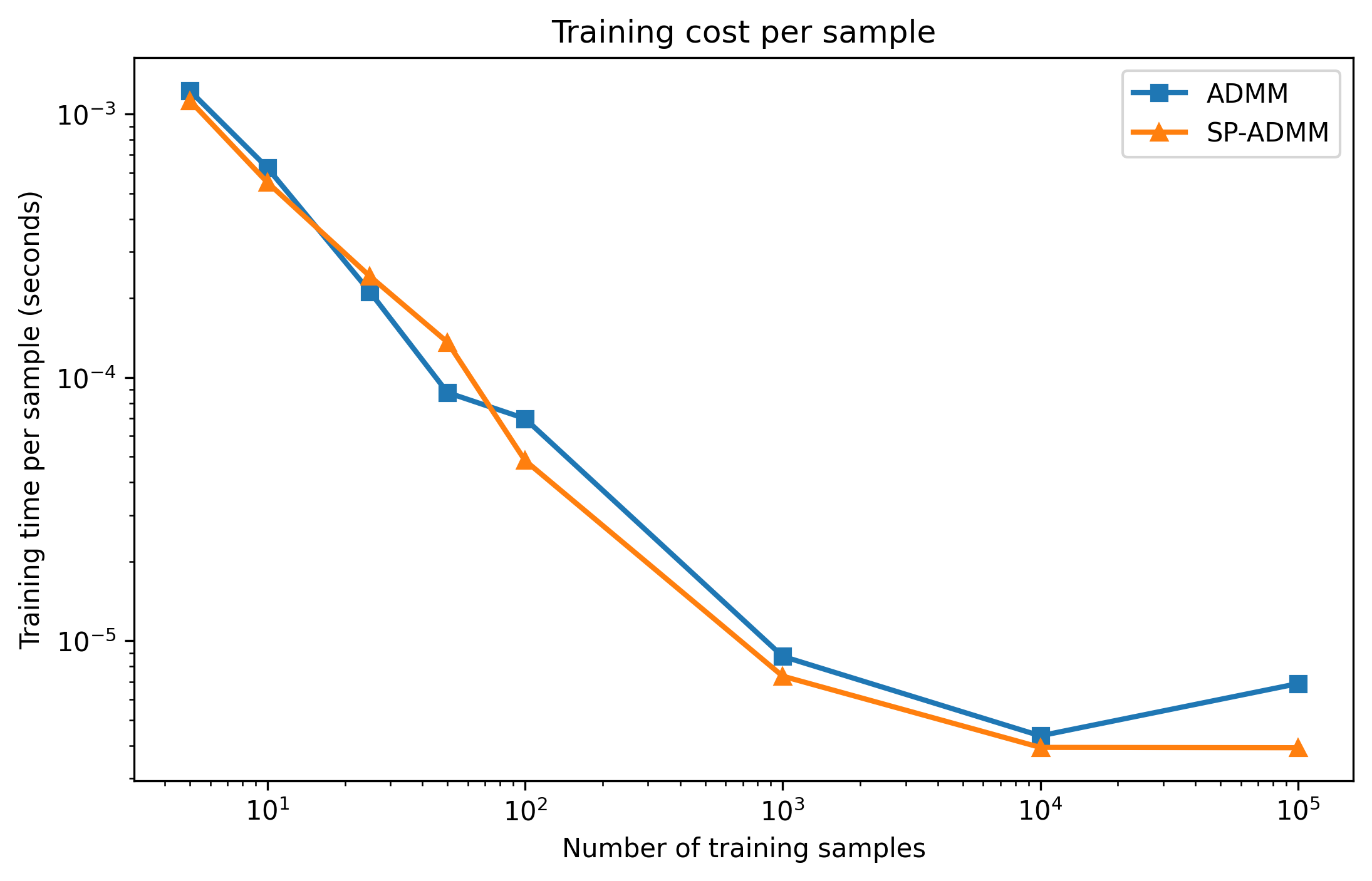}
    \caption{Training cost per sample.}
    \label{fig:training-time-b}
\end{subfigure}
\hfill
\begin{subfigure}[t]{0.32\textwidth}
    \centering
    \includegraphics[width=\textwidth]{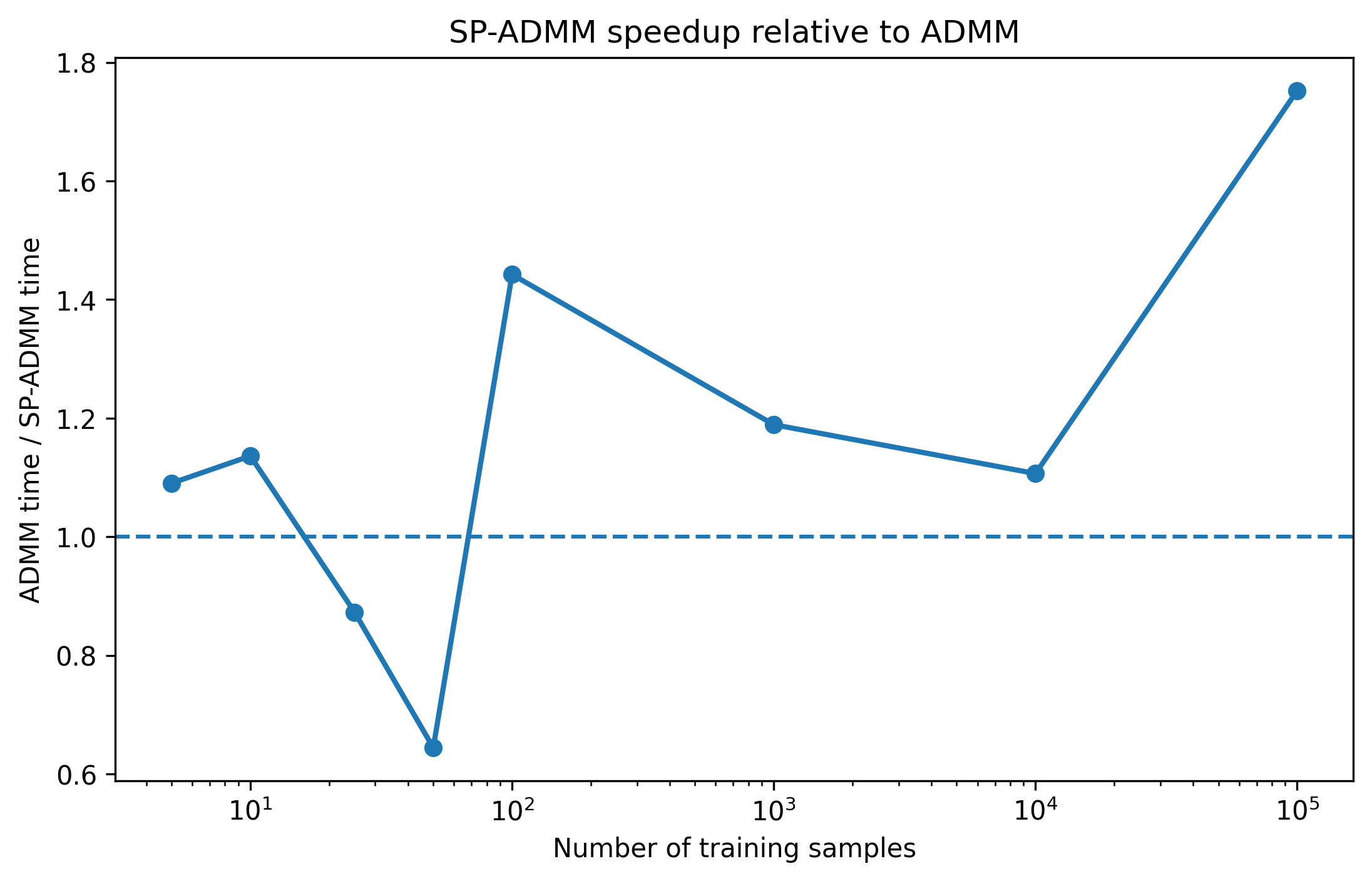}
    \caption{SP-ADMM speedup relative to ADMM.}
    \label{fig:training-time-c}
\end{subfigure}

\caption{Training-time comparison between ADMM and SP-ADMM. 
(a) Total wall-clock training time versus number of training samples. 
(b) Training time normalized per sample. 
(c) Relative speedup, measured as ADMM time divided by SP-ADMM time; values
above one indicate that SP-ADMM is faster.}
\label{fig:training-time-all}
\end{figure}

We also compare the computational cost of ADMM and SP-ADMM. The comparison
uses the same training data, noise level, stencil radius, and stopping
criteria as the training-size experiment. Since the finite-difference stencil
is fixed, FD has no training cost and is not included in the timing comparison.
Figure~\ref{fig:training-time-all} reports three timing diagnostics: total
wall-clock training time, training cost per sample, and the relative speedup
of SP-ADMM compared with ADMM.

Figure~\ref{fig:training-time-a} shows that the total training time increases
with the number of training samples for both learned methods. This is expected
because larger training sets produce larger least-squares systems. For small
training sets, the two methods have comparable runtime and may exchange order
because fixed overhead and solver setup costs are significant. For larger
training sets, SP-ADMM becomes faster than ADMM. This is consistent with the
structure-preserving formulation: SP-ADMM parameterizes the stencil directly
inside the skew-adjoint class, reducing the number of free unknowns and
avoiding the need to enforce skew-adjointness as an external equality
constraint.

Figure~\ref{fig:training-time-b} gives a normalized view of the training cost.
The cost per sample decreases as the training set grows, showing that fixed
setup costs are amortized over larger datasets. For the largest training
sets, SP-ADMM has a smaller per-sample cost than ADMM.

Figure~\ref{fig:training-time-c} summarizes the relative speedup. Values above
one indicate that SP-ADMM is faster than ADMM. The speedup is not monotone for
small sample sizes because timing measurements are affected by solver
overhead, initialization cost, and small problem-size effects. However, for
the largest training set, SP-ADMM is substantially faster. Together with the
accuracy results in Figure~\ref{fig:training-size}, these timing experiments
show that SP-ADMM is not only more accurate in many regimes, but also
computationally competitive and increasingly efficient for large training
sets.
\subsection{Effect of regularization}
\label{subsec:regularization}

\begin{figure}[htbp]
\centering

\begin{subfigure}[t]{0.48\textwidth}
    \centering
    \includegraphics[width=\textwidth]{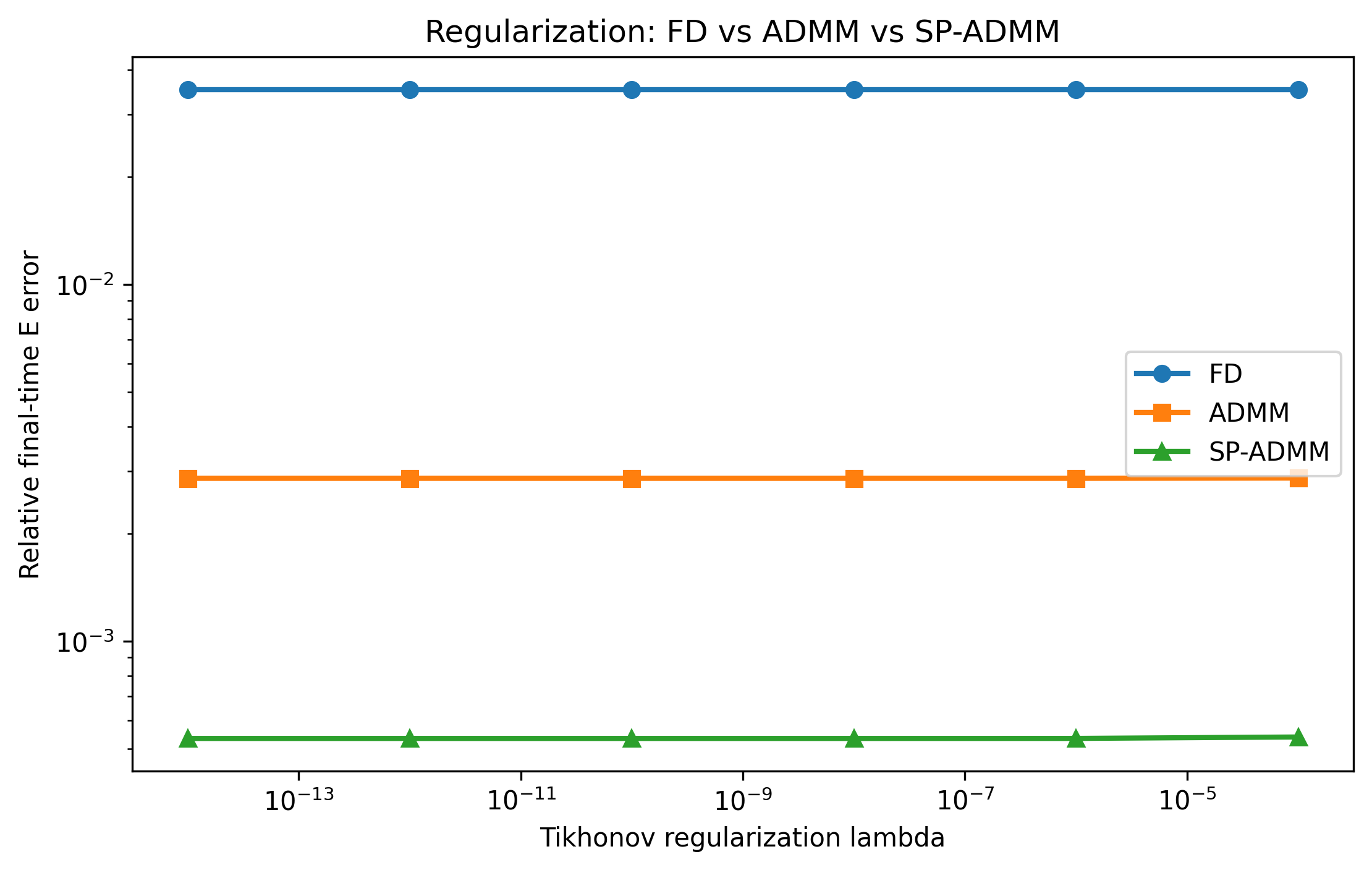}
    \caption{Final-time field error.}
    \label{fig:reg-a}
\end{subfigure}
\hfill
\begin{subfigure}[t]{0.48\textwidth}
    \centering
    \includegraphics[width=\textwidth]{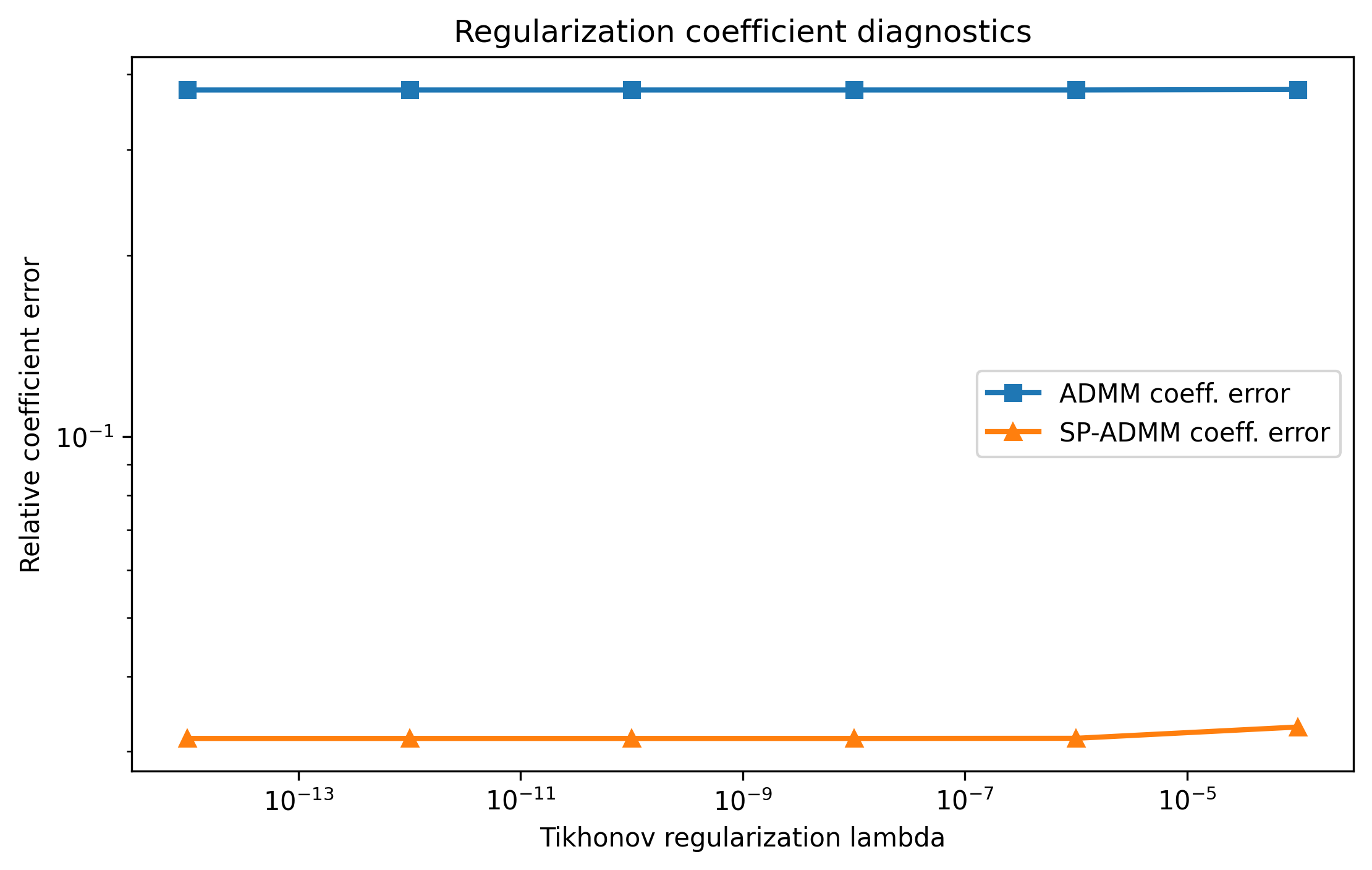}
    \caption{Coefficient error.}
    \label{fig:reg-b}
\end{subfigure}

\caption{Effect of Tikhonov regularization on learned-stencil performance.
(a) Final-time electric-field error for FD, ADMM, and SP-ADMM.
(b) Coefficient-recovery error for ADMM and SP-ADMM.
SP-ADMM gives the smallest field and coefficient errors over a broad range of
small-to-moderate regularization values, while overly large regularization
degrades both learned methods.}
\label{fig:regularization-all}
\end{figure}

The Tikhonov parameter $\lambda$ affects both field error and coefficient
recovery. Figure~\ref{fig:regularization-all} summarizes the effect of
regularization on the learned stencils: Figure~\ref{fig:reg-a} shows the
final-time field error, while Figure~\ref{fig:reg-b} shows the corresponding
coefficient-recovery behavior.

Figure~\ref{fig:reg-a} shows that the fixed FD error is unchanged with
respect to $\lambda$, while producing the largest error, since FD is not learned from data. In contrast, the
ADMM and SP-ADMM errors depend on the amount of regularization. Over a broad
range of small and moderate values of $\lambda$, SP-ADMM gives the smallest
final-time field error, while ADMM also improves substantially over FD.
However, when $\lambda$ becomes too large, the regularization begins to
over-penalize the learned coefficients, and the field error increases for
both learned methods.

Figure~\ref{fig:reg-b} shows a similar trend for coefficient recovery.
SP-ADMM has a noticeably smaller coefficient error than ADMM across most of
the regularization range, especially for small and moderate $\lambda$. This
suggests that the skew-parameterized formulation uses regularization more
effectively. For very large $\lambda$, both methods degrade, confirming that
excessive regularization can bias the learned stencil too strongly. Overall,
these results show that a moderate Tikhonov parameter improves robustness,
and that SP-ADMM is consistently less sensitive to the choice of
regularization than standard ADMM.
\subsection{Constraint ablation}
\label{subsec:ablation}

\begin{figure}[htbp]
\centering

\begin{subfigure}[t]{0.48\textwidth}
    \centering
    \includegraphics[width=\textwidth]{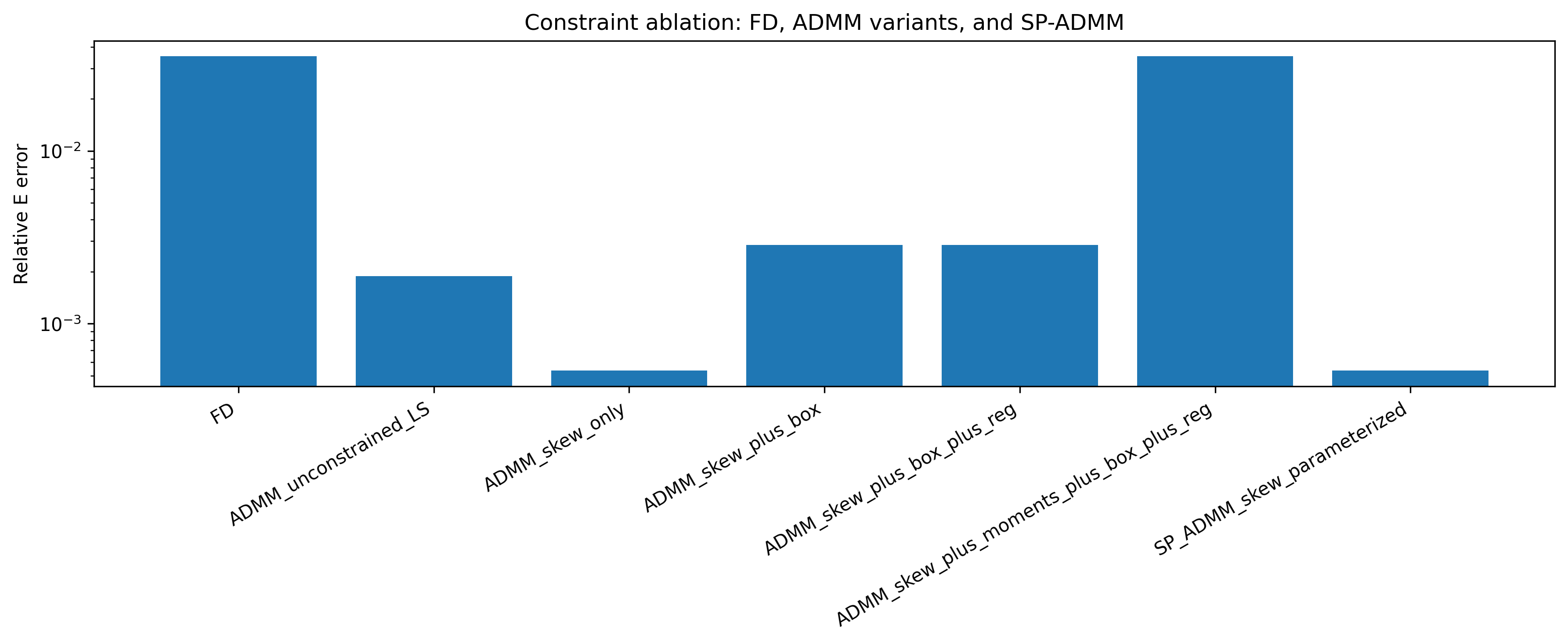}
    \caption{Field error.}
    \label{fig:ablation-a}
\end{subfigure}
\hfill
\begin{subfigure}[t]{0.48\textwidth}
    \centering
    \includegraphics[width=\textwidth]{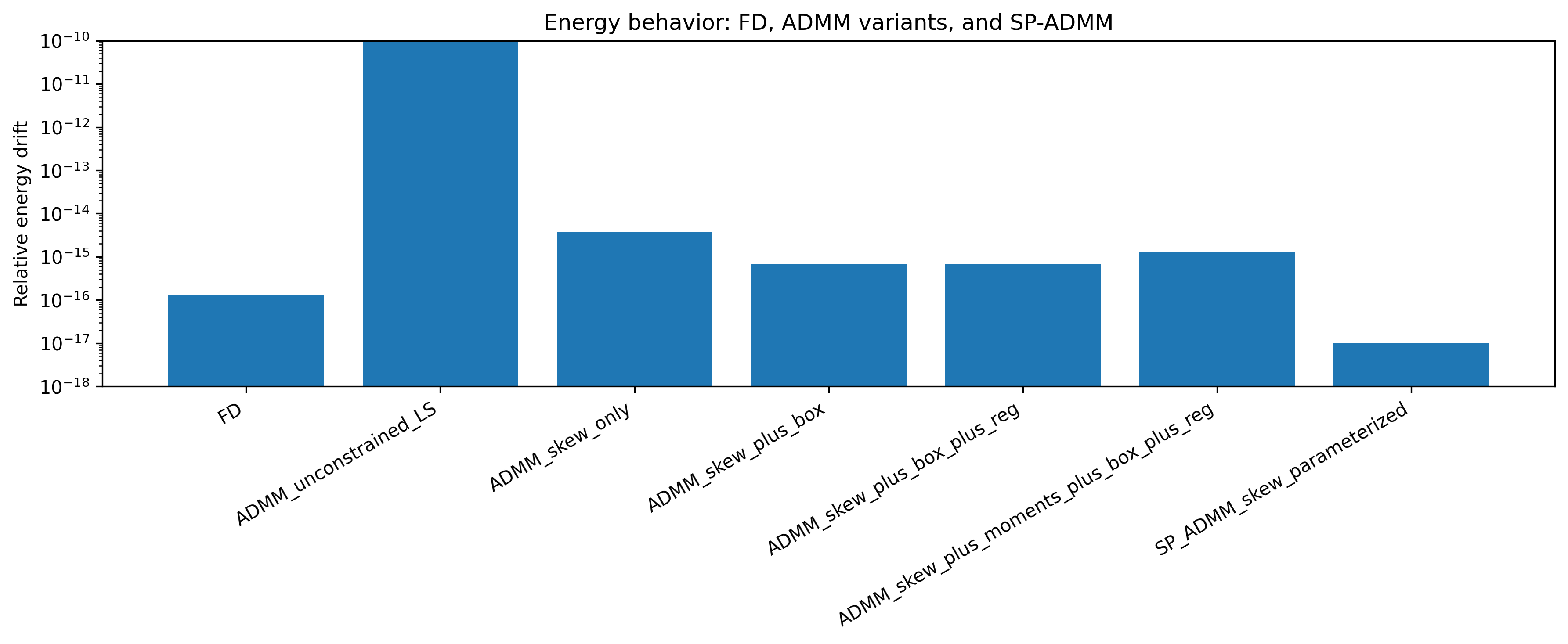}
    \caption{Energy drift.}
    \label{fig:ablation-b}
\end{subfigure}

\caption{Constraint ablation for FD, unconstrained least squares, ADMM
variants, and SP-ADMM. (a) Final-time field error. The skew-only ADMM and
SP-ADMM variants give the smallest errors, while imposing moment constraints
can increase the error in the hidden-operator regime by forcing the stencil
toward classical finite-difference behavior. (b) Energy behavior. The
unconstrained least-squares fit shows much larger energy drift, whereas the
skew-constrained and skew-parameterized methods preserve energy to roundoff
accuracy.}
\label{fig:ablation-all}
\end{figure}

The ablation experiment compares FD, unconstrained least squares, several
ADMM variants, and SP-ADMM. Figure~\ref{fig:ablation-all} summarizes the
results: Figure~\ref{fig:ablation-a} compares the final-time field errors,
while Figure~\ref{fig:ablation-b} shows the corresponding energy drift.

Figure~\ref{fig:ablation-a} shows that skew-adjointness is the most important
constraint for accuracy in the hidden-operator regime. The unconstrained
least-squares fit improves over FD, but the best results are obtained from
the skew-only ADMM and SP-ADMM formulations. By contrast, adding moment
constraints increases the error in this experiment, because those constraints
bias the learned stencil toward classical finite-difference coefficients
rather than allowing it to adapt to the hidden nonstandard skew-adjoint
operator.

Figure~\ref{fig:ablation-b} shows that unconstrained least squares also has
substantially worse energy behavior. In contrast, all skew-constrained and
skew-parameterized methods preserve the discrete Maxwell energy to roundoff
accuracy. This confirms that skew-adjointness is the essential structure for
energy conservation, while moment constraints from taylor construction are mainly useful for
finite-difference recovery rather than for learning a hidden nonstandard
operator.

\subsection{Error growth over increasing final times}
\label{subsec:time}

\begin{figure}[htbp]
\centering

\begin{subfigure}[t]{0.48\textwidth}
    \centering
    \includegraphics[width=\textwidth]{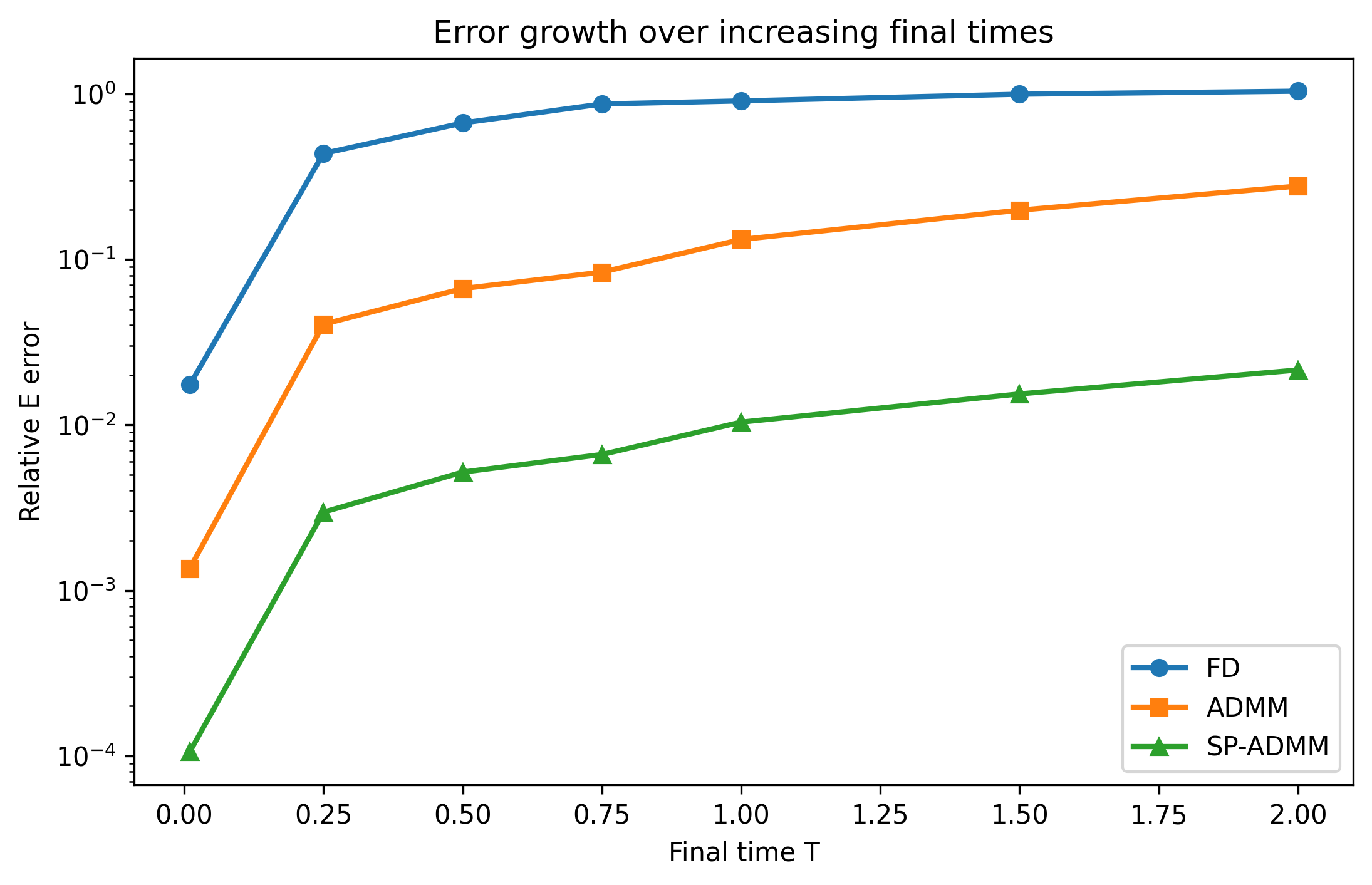}
    \caption{Final-time field error.}
    \label{fig:time-a}
\end{subfigure}
\hfill
\begin{subfigure}[t]{0.48\textwidth}
    \centering
    \includegraphics[width=\textwidth]{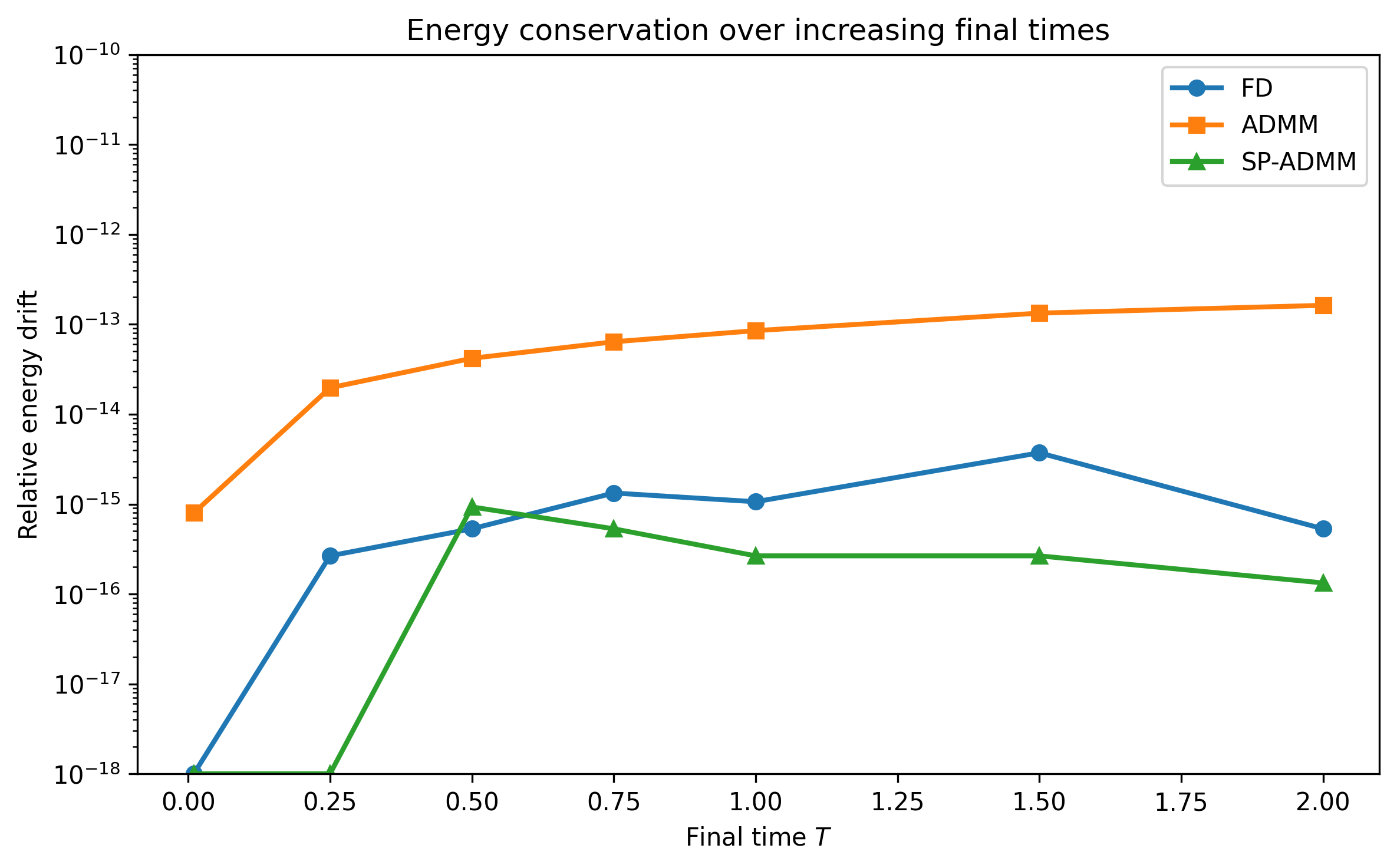}
    \caption{Energy drift.}
    \label{fig:time-b}
\end{subfigure}

\caption{Long-time behavior as the final time increases. (a) Final-time
electric-field error. FD exhibits the fastest error growth, ADMM performs
better than FD, and SP-ADMM remains the most accurate over the full time
interval. (b) Energy conservation over increasing final times. All methods
preserve the discrete Maxwell energy to near roundoff accuracy; values below
$10^{-17}$ are displayed at the floor for log-scale visualization.}
\label{fig:time-all}
\end{figure}

Finally, we test whether the error advantage persists as the final time
increases. Figure~\ref{fig:time-all} summarizes the long-time behavior:
Figure~\ref{fig:time-a} shows the growth of the final-time electric-field
error, while Figure~\ref{fig:time-b} shows the corresponding energy drift.

Figure~\ref{fig:time-a} shows that the advantage of the learned stencils
persists as the final time increases. The FD error grows the fastest and
reaches order one over the tested interval. ADMM reduces this growth
substantially, while SP-ADMM remains the most accurate method at all final
times shown. This indicates that the structure-preserving learned stencil
produces a more faithful long-time evolution than both the fixed
finite-difference stencil and the standard ADMM-learned stencil.

Figure~\ref{fig:time-b} shows that all three methods maintain very small
energy drift, remaining close to roundoff level throughout the simulation.
SP-ADMM exhibits the smallest overall drift for most final times, while ADMM
has a slightly larger but still negligible drift. This further highlights the long-time accuracy advantage of SP-ADMM.

\subsection{Layered-medium Maxwell wave propagation}
\label{subsec:layered-medium}

To complement the hidden-operator experiments, we consider a physically
motivated Maxwell propagation problem in a one-dimensional layered medium.
This example is included to test the learned structure-preserving stencils
in a setting that produces reflection and transmission at a material
interface, rather than only in periodic hidden-operator learning tests.

We consider the variable-coefficient Maxwell system
\begin{equation}
\begin{aligned}
  \varepsilon(x) E_t &= H_x,\\
  \mu H_t &= E_x,
\end{aligned}
\qquad x\in[0,L], \quad t>0,
\label{eq:layered-maxwell}
\end{equation}
where $\mu=1$ and the permittivity is piecewise constant,
\begin{equation}
  \varepsilon(x)=
  \begin{cases}
  \varepsilon_1, & 0\le x < x_I,\\
  \varepsilon_2, & x_I\le x \le L.
  \end{cases}
  \label{eq:layered-eps}
\end{equation}
Here $x_I$ denotes the material interface. The jump in $\varepsilon(x)$ changes
the wave speed and impedance, causing an incident pulse to separate into
reflected and transmitted waves. Such layered-medium tests are standard
benchmarks in computational electromagnetics and wave propagation because they
probe the behavior of a numerical method near material interfaces
\cite{yee1966numerical,taflove2000computational,leveque2007finite}.

The initial electric field is chosen as a localized Gaussian pulse in the left
material,
\begin{equation}
  E(x,0)=\exp\left(-\frac{(x-x_0)^2}{\sigma^2}\right),
  \label{eq:layered-initial-E}
\end{equation}
and the magnetic field is initialized consistently as a right-going wave in
the first material,
\begin{equation}
  H(x,0)=-\sqrt{\frac{\varepsilon_1}{\mu}}\,E(x,0).
  \label{eq:layered-initial-H}
\end{equation}
 We compare FD, ADMM, and SP-ADMM against a
high-resolution reference solution. The learned stencils are trained using a
noise level $\delta=0.05$ and Tikhonov regularization
$\lambda=10^{-8}$.

\begin{figure}[htbp]
\centering

\begin{subfigure}[t]{0.95\textwidth}
    \centering
    \includegraphics[width=\textwidth]{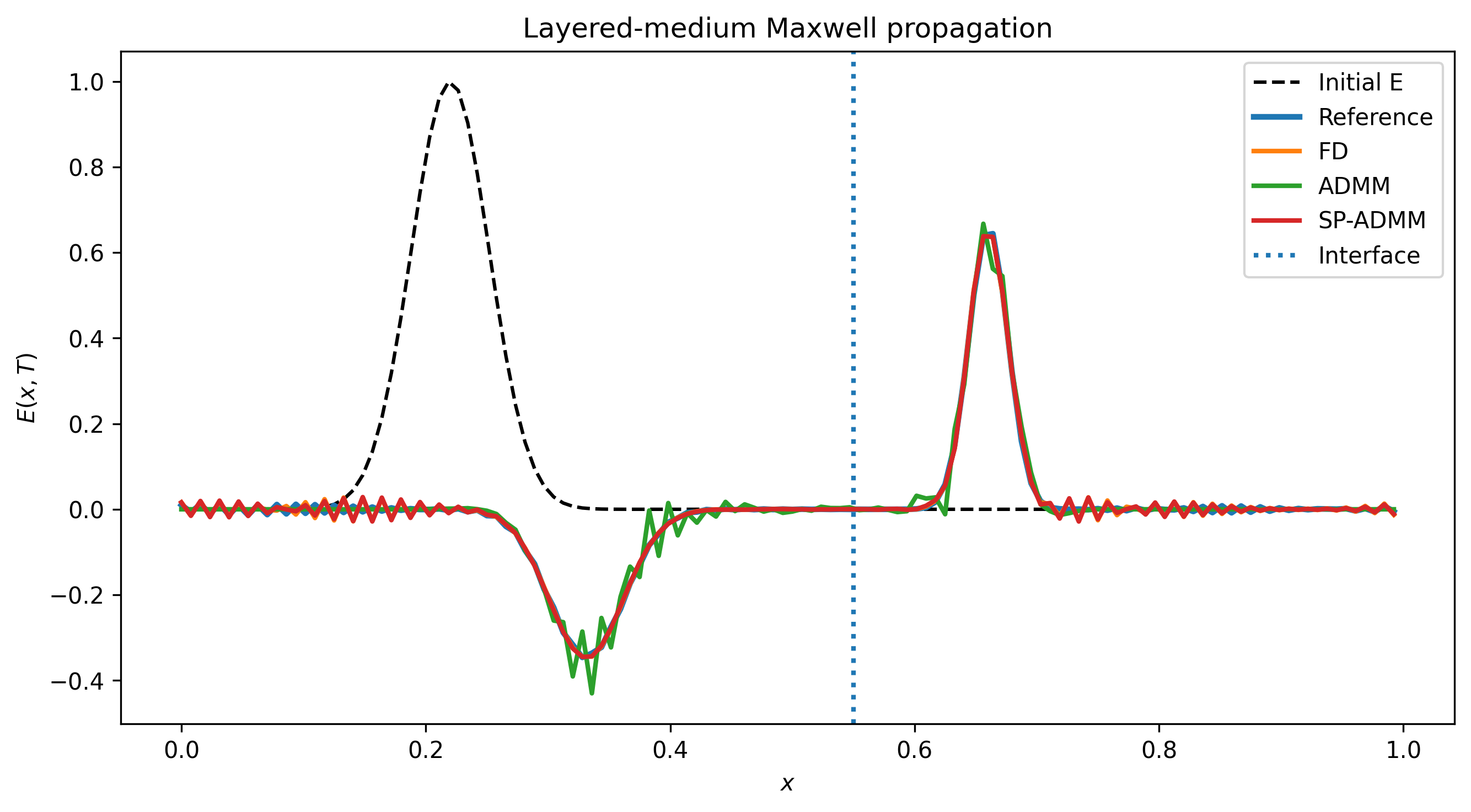}
    \caption{Electric-field profile at final time.}
    \label{fig:layered-profile}
\end{subfigure}

\vspace{0.25cm}

\begin{subfigure}[t]{0.48\textwidth}
    \centering
    \includegraphics[width=\textwidth]{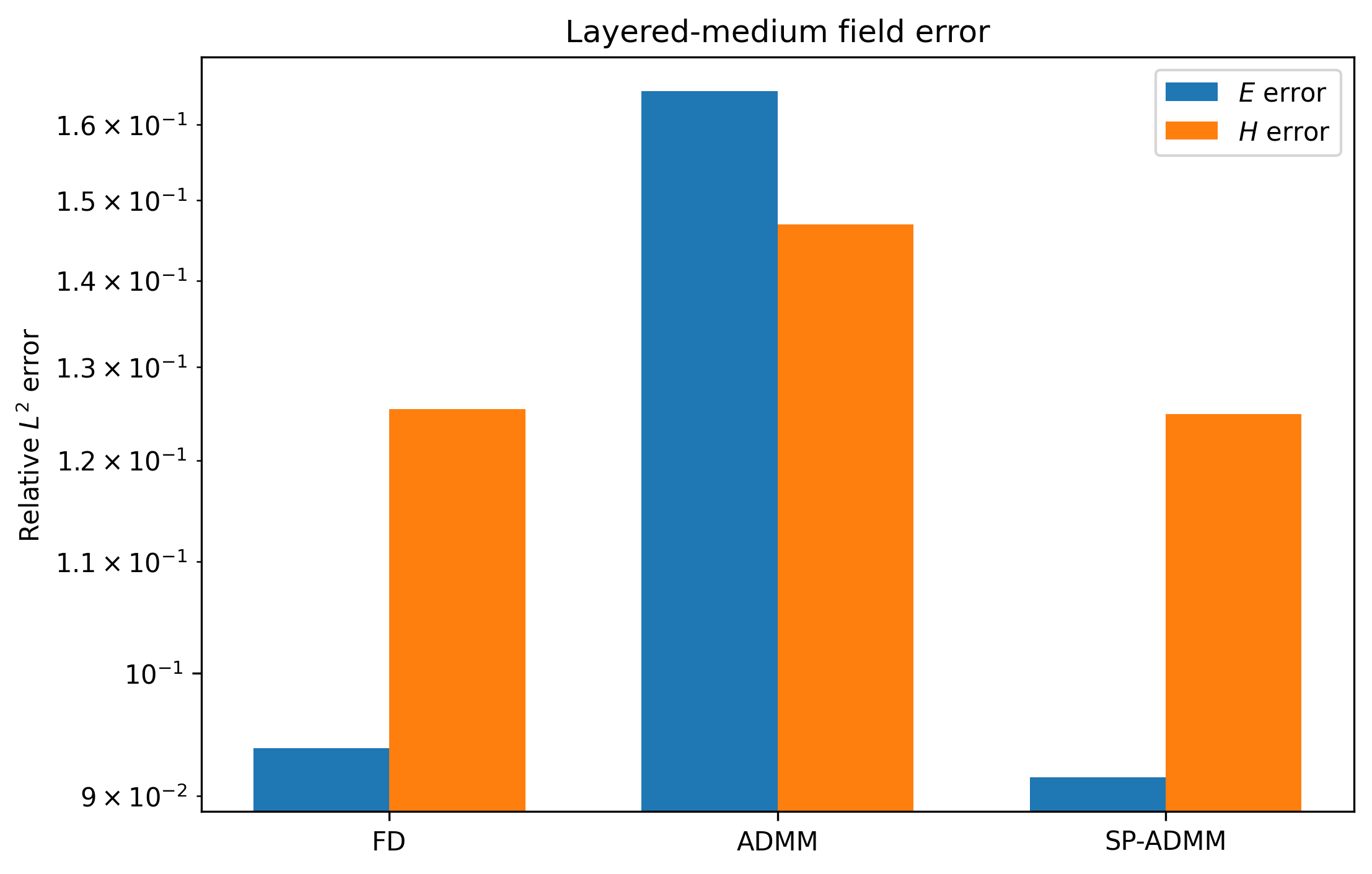}
    \caption{Relative $L^2$ field errors.}
    \label{fig:layered-error}
\end{subfigure}
\hfill
\begin{subfigure}[t]{0.48\textwidth}
    \centering
    \includegraphics[width=\textwidth]{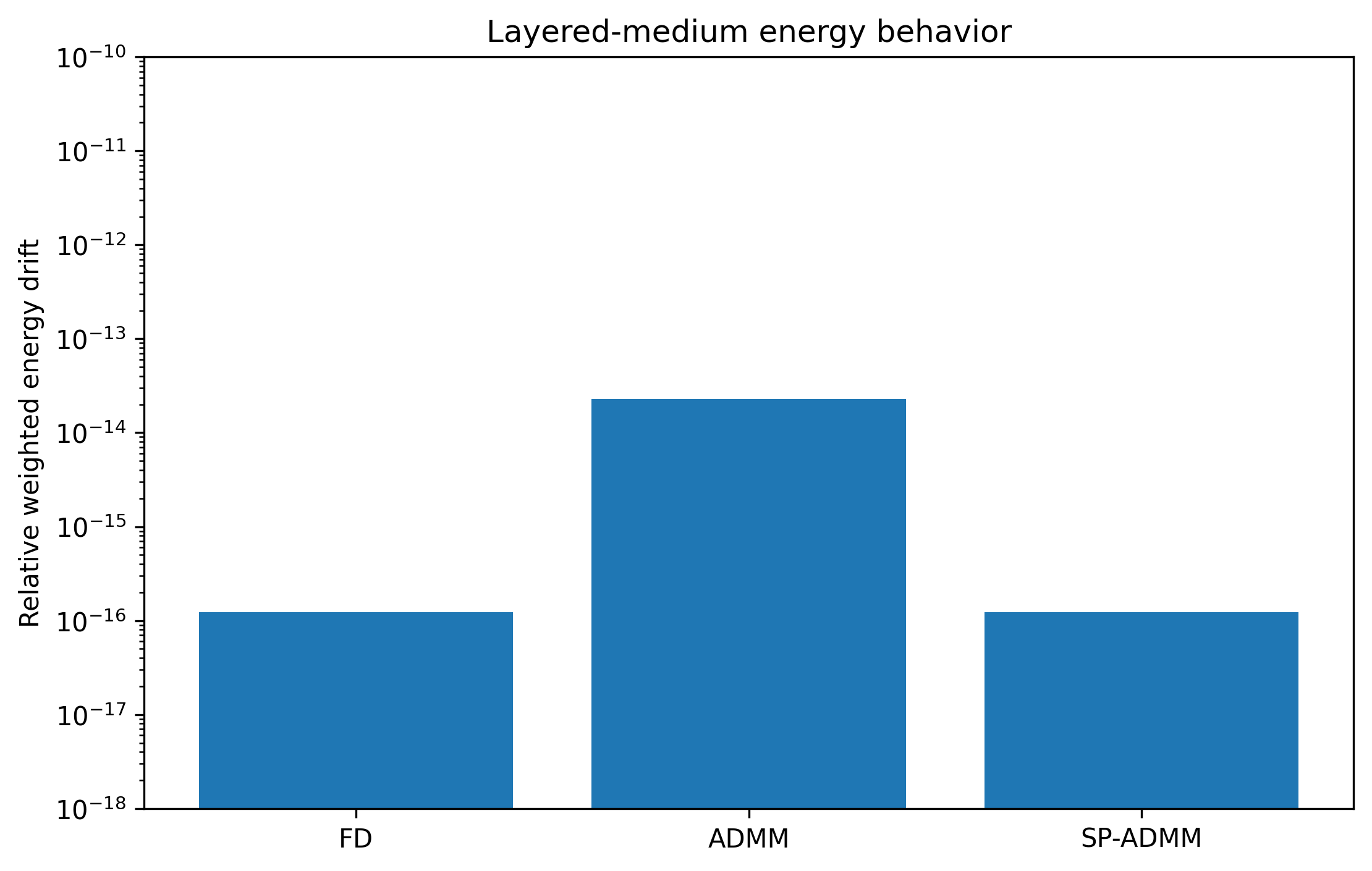}
    \caption{Relative weighted energy drift.}
    \label{fig:layered-energy-plot}
\end{subfigure}

\caption{Layered-medium Maxwell propagation test. A Gaussian pulse travels
toward a dielectric interface and separates into reflected and transmitted
waves. FD and SP-ADMM both closely follow the reference solution, while ADMM
shows larger oscillations near the wave fronts. The error plot shows that FD
and SP-ADMM are quantitatively comparable, while ADMM gives the largest field
errors. The energy plot shows that FD and SP-ADMM preserve the weighted
electromagnetic energy to near roundoff accuracy.}
\label{fig:layered-medium-all}
\end{figure}

Figure~\ref{fig:layered-medium-all} shows the layered-medium Maxwell
propagation test. The profile plot shows that the incident Gaussian pulse has
reached the material interface and split into reflected and transmitted wave
components. All three methods reproduce the main wave structure. The FD and
SP-ADMM curves remain close to the reference solution, while ADMM exhibits
more visible oscillations near the reflected and transmitted wave fronts.
 
ADMM has the largest error in both fields. Thus, in this
standard layered-medium test, the learned SP-ADMM stencil remains competitive
with the classical finite-difference stencil.
The energy plot shows that FD and SP-ADMM preserve the
electromagnetic energy to near roundoff accuracy, while ADMM has a larger,
but still small, energy drift. This confirms that the small oscillations
visible in the SP-ADMM solution should not be interpreted as energy
instability. The SP-ADMM stencil remains skew-adjoint and energy-preserving,
but skew-adjointness alone does not imply a non-oscillatory approximation.

The oscillations near the material interface are consistent with the known
behavior of nondissipative high-order or spectral-type approximations near
sharp transitions. In the layered-medium problem, the discontinuity in
$\varepsilon(x)$ reduces the local smoothness of the solution. Gibbs-type
oscillations near discontinuities are well documented for high-order
approximations, and common remedies include limiters, filtering, or artificial
dissipation \cite{gottlieb1997gibbs,leveque2007finite,discacciati2020viscosity}.
Therefore, this experiment highlights both the strength and limitation of the
proposed approach: SP-ADMM preserves the energy structure and remains
competitive with FD in a physical Maxwell propagation problem, but additional
interface-aware stabilization may be useful when sharp material jumps are
present.
\section{Conclusion and Future Work}
\label{sec:conclusion}

This paper studied structure-preserving stencil learning for the
one-dimensional periodic Maxwell system. We compared a fixed finite-difference
(FD) stencil, a constrained ADMM-learned stencil, and a structure-preserving
ADMM method (SP-ADMM) that enforces skew-adjointness directly through the
stencil parameterization.

The results show that learned stencils can outperform fixed finite
differences when the data are generated by a nonstandard hidden
skew-adjoint operator. Across clean data, noisy derivative data, different
initial conditions, and several hidden operators, SP-ADMM gave the
smallest final-time electric-field error. At the same time, it preserves the
discrete Maxwell energy to roundoff accuracy.

The ablation study shows that skew-adjointness is the essential constraint
for energy conservation.Moment constraints are useful for finite-difference recovery, but they can reduce accuracy when the hidden operator is not a classical finite-difference stencil, because they bias the learned stencil toward standard finite-difference coefficients. The training-size, regularization, and timing tests also
show that SP-ADMM benefits from additional data and remains computationally
competitive with ADMM.

Overall, SP-ADMM provides a data-driven way to learn accurate Maxwell
stencils while preserving the energy-conserving structure of the equations.
Future work will extend this approach to higher-dimensional Maxwell systems,
higher-order stencils, nonlinear dispersive media, and additional structure
constraints such as discrete divergence preservation.

\section*{Availability of Data and Code}
The code used to generate the numerical results in this study is
available at
\url{https://github.com/victoryobieke/energy-stable-maxwell-stencil-learning}.
\section*{Conflict of Interest}
The authors declare that there are no conflicts of interest.

\section*{Funding}
Not applicable. The author received no specific funding for this work.

\section*{Author Contributions}
The authors are solely responsible for the conceptualization, methodology, software implementation, numerical experiments, analysis, and writing of this manuscript.
\bibliographystyle{plain}
\bibliography{refs}

\begin{thebibliography}{99}

\bibitem{yee1966}
K.~S. Yee.
\newblock Numerical solution of initial boundary value problems involving Maxwell's equations in isotropic media.
\newblock \emph{IEEE Transactions on Antennas and Propagation}, 14(3):302--307, 1966.

\bibitem{taflove2005}
A.~Taflove and S.~C. Hagness.
\newblock \emph{Computational Electrodynamics: The Finite-Difference Time-Domain Method}.
\newblock Artech House, 3rd edition, 2005.
\bibitem{obieke2025structure}
V.~Obieke and E.~Oguadimma.
\newblock Structure-preserving physics-informed neural network for the
Korteweg--de Vries (KdV) equation.
\newblock \emph{arXiv preprint arXiv:2511.00418}, 2025.
\bibitem{hairer2006}
E.~Hairer, C.~Lubich, and G.~Wanner.
\newblock \emph{Geometric Numerical Integration: Structure-Preserving Algorithms for Ordinary Differential Equations}.
\newblock Springer, 2nd edition, 2006.

\bibitem{boyd2011admm}
S.~Boyd, N.~Parikh, E.~Chu, B.~Peleato, and J.~Eckstein.
\newblock Distributed optimization and statistical learning via the alternating direction method of multipliers.
\newblock \emph{Foundations and Trends in Machine Learning}, 3(1):1--122, 2011.

\bibitem{barsinai2019}
Y.~Bar-Sinai, S.~Hoyer, J.~Hickey, and M.~P. Brenner.
\newblock Learning data-driven discretizations for partial differential equations.
\newblock \emph{Proceedings of the National Academy of Sciences}, 116(31):15344--15349, 2019.

\bibitem{brunton2016}
S.~L. Brunton, J.~L. Proctor, and J.~N. Kutz.
\newblock Discovering governing equations from data by sparse identification of nonlinear dynamical systems.
\newblock \emph{Proceedings of the National Academy of Sciences}, 113(15):3932--3937, 2016.
\bibitem{bokil2019highorder}
V.~A. Bokil and N.~L. Gibson.
\newblock High spatial order energy stable FDTD methods for Maxwell's equations
in nonlinear optical media in one dimension.
\newblock \emph{Journal of Scientific Computing}, 77(1):330--371, 2018.

\bibitem{obieke2025structure}
V.~Obieke and E.~Oguadimma.
\newblock Structure-preserving physics-informed neural network for the
Korteweg--de Vries (KdV) equation.
\newblock \emph{arXiv preprint arXiv:2511.00418}, 2025.
\bibitem{higham2008functions}
N. J. Higham,
\textit{Functions of Matrices: Theory and Computation},
SIAM, Philadelphia, 2008.

\bibitem{moler2003nineteen}
C. Moler and C. Van Loan,
Nineteen dubious ways to compute the exponential of a matrix,
twenty-five years later,
\textit{SIAM Review}, 45(1), pp. 3--49, 2003.

\bibitem{gottlieb1997gibbs}
D.~Gottlieb and C.-W.~Shu.
\newblock On the Gibbs phenomenon and its resolution.
\newblock \emph{SIAM Review}, 39(4):644--668, 1997.

\bibitem{leveque2002finite}
R.~J.~LeVeque.
\newblock \emph{Finite Volume Methods for Hyperbolic Problems}.
\newblock Cambridge University Press, Cambridge, 2002.

\bibitem{discacciati2020oscillations}
N.~Discacciati, J.~S.~Hesthaven, and D.~Ray.
\newblock Controlling oscillations in high-order discontinuous Galerkin schemes
using artificial viscosity tuned by neural networks.
\newblock \emph{Journal of Computational Physics}, 409:109304, 2020.

\end{thebibliography}

\end{document}